\documentclass[12pt]{article}
\usepackage[english]{babel}
\usepackage{graphics}
\usepackage[all]{xy}
\usepackage{epsfig}
\usepackage{amsmath,amssymb,latexsym}
\usepackage{color,graphicx,indentfirst}
\def\a{\alpha}  \def\b{\beta}        
\def\de{\delta}  
 \def\vk{\varkappa}

\def\bPhi{\boldsymbol\Phi}
\def\ov{\overline}
\def\wh{\widehat}
\def\wt{\widetilde}
\def\ds{\displaystyle}
\def\pa{\partial}
\def\thalf{\textstyle{\frac{1}{2}}\displaystyle}

\def\ov{\overline}
\def\pa{\partial}

\def\un{\underline}
\def\wh{\widehat}
\def\wt{\widetilde}
\def\ds{\displaystyle}

\def\lb{\label}
\def\nn{\nonumber}
\def\lf{\left}
\def\rg{\right}
\def\leq{\leqslant}
\def\geq{\geqslant}
\DeclareMathOperator{\Rea}{Re}
\DeclareMathOperator{\Ima}{Im}
\DeclareMathOperator{\erfc}{erfc}
\newtheorem{proposition}{\indent Proposition}
\newtheorem{lemma}{\indent Lemma}
\newtheorem{corollary}{\indent Corollary}
\newtheorem{remark}{\indent Remark}
\begin{document}
\large
\centerline{\bf On a family of finite-difference schemes with discrete}
\smallskip
\centerline{\bf transparent boundary conditions}
\smallskip
\centerline{\bf for a parabolic equation on the half-axis}
\smallskip
\begin{center}{Alexander Zlotnik and Natalya Koltsova}
\end{center}
\smallskip
\begin{center}{\it Department of Mathematics at Faculty of Economics,
National Research University Higher School of Economics,
Myasnitskaya 20, 101990 Moscow, {\sc Russia}}
\end{center}
\smallskip
\begin{center}{\it Department of Mathematical Modelling,
National Research University Moscow Power Engineering Institute,
Krasnokazarmennaya 14, 111250 Moscow, {\sc Russia}}
\end{center}
\smallskip
\centerline{\it e-mail: azlotnik2008@gmail.com, predienloko@yandex.ru}
\begin{abstract}
\noindent An initial-boundary value problem for the 1D self-adjoint parabolic equation on the half-axis is solved.
We study a broad family of two-level finite-difference schemes with two parameters
related to averagings both in time and space.
Stability in two norms is proved by the energy method.
Also discrete transparent boundary conditions are rigorously derived for schemes by applying the method of reproducing functions.
Results of numerical experiments are included as well.
\end{abstract}
\bigskip
\section{Introduction}
\label{sect1}
In many applications, a problem of solving partial differential equations
in unbounded domains arises.
A number of approaches to the problem is developed mainly associated with the statement of additional boundary conditions on artificial boundaries \cite{G87}-\cite {ME05}.
The conditions are called the (exact) artificial/non-reflecting/transparent boundary conditions provided that they are satisfied by the solutions of the original problems in unbounded domains. For definiteness, we exploit the last name (TBCs).
For parabolic evolution equations or for the Schro\"{o}dinger equation, the TBCs are integro-differential relations along the artificial boundaries.
Their adequate discretization is non-trivial since it can produce significant reflections from the boundaries and even instability in computations as well as create difficulties in rigorous proofs of stability of the resulting numerical method.
\par An alternative approach suggests to implement the idea of the TBCs on the mesh level. Namely, first one consider a discretization of the problem on an infinite mesh in the unbounded domain (which is not practical because of the infinite number of unknowns). Its solution is restricted to the finite mesh by deriving a mesh counterpart of the TBC on the artificial boundary.
One version of this approach is associated with the derivation of \textit{discrete} TBCs requiring to solve analytically model mesh problems on infinite grids.
This approach had worked well by the complete absence of reflections from artificial boundaries and reliable stability of computations in practice as well as by clarity of the mathematical background and rigorous proofs of stability of the resulting mesh method in theory.
Such an approach is developed in detail for 1D time-dependent Schr\"{o}dinger equation in \cite {EA01}-\cite {DZZ09} and also used for 1D parabolic equations in \cite{G87,ME0,ME01,ME05,Z07}.
\par In this paper, the approach is developed for the 1D self-adjoint parabolic equation on the half-axis. We study a broad family of two-level finite-difference $\sigma$-schemes with averaging in space with a weight $\theta$.
We prove its stability by the energy method and rigorously derive the discrete TBC by applying the method of reproducing functions. Notice that both points are not so well developed in some other papers.
Similarly to \cite{DZZ09}, this allows one to cover in a unified manner a collection of particular schemes: the standard scheme without averaging ($\theta=0$),
the linear finite-element method ($\theta=\frac16$),
scheme of higher order of accuracy (for constant coefficients, $\theta=\frac{1}{12}$)
and a vector scheme on a four-point stencil ($\theta=\frac14$).
The results generalize those obtained for $\theta=0$ in \cite{Z07} but
here we prove stability in two (not one) energy norms, the rigorous derivation of the discrete TBC is notably different and results of numerical experiments
are included as well.
\par Notice that the results for the family of finite-difference schemes can be enlarged for the 2D (or multi-D) case exploiting the technique from \cite{ZZ11,IZ11}.
\section{An initial-boundary value problem and
a finite-difference $\sigma$-scheme with averaging in space and an approximate TBC}
\label{sect2}
We consider the one-dimensional parabolic equation
\begin{equation}
 \rho\,\frac{\pa u}{\pa t}+{\mathcal A}u=f,\ \
 {\mathcal A}u:=-\frac{\pa}{\pa x}
 \lf(b\frac{\pa u}{\pa x}\rg)+cu
\label{sch}
\end{equation}
for $x>0$ and $t>0$. Its coefficients satisfy the conditions
$\rho(x)\geq \un{\rho}>0$, $b(x)\geq \nu>0$ and $c(x)\geq 0$ for $x>0$.
\par Equation (\ref{sch}) is supplemented with the following boundary condition, the condition at infinity and the initial condition:
\begin{gather}
u|_{x=0}=g(t),\ \ u(x,t)\to 0\ \ \mbox{as}\ x\to \infty, \ \mbox{for all}\ t>0,
 \label{bc}\\[1mm]
\lf. u\rg|_{t=0}=u^0(x)\ \ \mbox{for}\ x>0.
\label{ic}
\end{gather}
We assume that the coefficients become constants as well as
$f$ and $u^0$ vanish for sufficiently large $x\geq X_0$ (for some $X_0>0$)
\begin{equation}
 \rho(x)=\rho_\infty>0,\ b(x)=b_\infty>0,\ c(x)=c_{\infty}\geq 0,\
 f(x,t)=0,\ u^0(x)=0.
\label{cc}
\end{equation}
\par An integro-differential TBC satisfied by the solution to this problem can be written in the Dirichlet-to-Neumann form
\begin{equation}
 \frac{\pa u}{\pa x}\,(X,t)
 =-\sqrt{\frac{\rho_\infty}{b_\infty}}\,e^{\ds{-(c_\infty/\rho_\infty)t}}\,
 \frac{1}{\sqrt{\pi}}\frac{d}{dt}\int_0^t
 u(X,\theta) e^{\ds{(c_\infty/\rho_\infty)\theta}}\,
 \frac{d\theta}{\sqrt{t-\theta}}
\label{tbc}
\end{equation}
for $t>0$ and for any $X\geq X_0$; other equivalent forms are also known.
The TBC is nonlocal in time; recall that the involved operator
\[
 \mathcal{D}_{t+}^{1/2}w(t):=\frac{1}{\sqrt{\pi}}\frac{d}{dt}\int_0^t
 w(\theta)\, \frac{d\theta}{\sqrt{t-\theta}}\ \ \mbox{for}\ \ t>0
\]
defines the classical left-hand Riemann-Liouville time derivative of order
$\frac12$ on the half-axis $[0,\infty)$. But we do not exploit this TBC explicitly below.
\par We fix some $X>X_0$, set $\Omega=(0,X)$ and introduce a nonuniform mesh
$\ov{\omega}_{h,\infty}$ in $x$ on $[0,\infty)$ with nodes $0=x_0<\dots <x_J=X<\dots$
and steps $h_j:=x_j-x_{j-1}$ such that $h_J\leq X-X_0$ and $h_j=h\equiv h_J$
for $j\geq J$.
We set $x_{j-1/2}:=\frac{x_{j-1}+x_j}{2}$ and $h_{j+1/2}:=\frac{h_j+h_{j+1}}{2}$.
Let $\omega_{h,\infty}:=\ov{\omega}_{h,\infty}\backslash\,\{0\}$,
$\ov{\omega}_h:=\{x_j\}_{j=0}^J$ and $\omega_h:=\{x_j\}_{j=1}^{J-1}$.
Let $H_0(\ov{\omega}_h)$ be the space of functions on mesh $\ov{\omega}_h$ that equal $0$ for $x_0=0$.
\par We define the backward, modified forward and central difference quotients in $x$
\[
 \ov{\pa}_xW_j:= \frac{W_j-W_{j-1}}{h_j},\ \
 \wh{\pa}_xW_j:= \frac{W_{j+1}-W_{j}}{h_{j+1/2}},\ \
 \overset{\circ}{\pa}_x W_j:= \frac{W_{j+1}-W_{j-1}}{2h_{j+1/2}},
\]
together with averaging operators in $x$
\begin{gather*}
 \wh{s}_xW_j
 := \frac{h_j}{2h_{j+1/2}} W_j +\frac{h_{j+1}}{2h_{j+1/2}}W_{j+1},
\\[1mm]
 s_\theta W_j:=\theta\frac{h_j}{h_{j+1/2}} W_{j-1}+(1-2\theta) W_j+\theta\frac{h_{j+1}}{h_{j+1/2}} W_{j+1}
\end{gather*}
and a more general mesh counterpart of multiplication by a mesh function  $\varkappa_h$
\[
 C_\theta[\varkappa_h] W_j
 :=\theta\frac{h_j}{h_{j+1/2}}\varkappa_{hj} W_{j-1}
 +(1-2\theta)(\wh{s}_x\varkappa_{hj}) W_j
 +\theta\frac{h_{j+1}}{h_{j+1/2}}\varkappa_{h_{j+1}} W_{j+1},
\]
see \cite{DZZ09}; clearly $C_\theta[1]=s_\theta$. Notice that
$s_\theta W_J=s_\theta^{-} W_J+s_\theta^{+}W_J$, where
\[
 s_\theta^{-} W_J:=\theta W_{J-1}+(\thalf-\theta)W_J, \ \
 s_\theta^{+} W_J:=(\thalf-\theta)W_J+\theta W_{J+1}.
\]
\par We also introduce a nonuniform mesh in
$t$ on $[0,\infty)$ with the nodes $0=t_0<\dots<t_m<\dots$ such that $t_m\to\infty$ for $m\to\infty$
and the steps $\tau_m:=t_m-t_{m-1}$.
Let
$\omega^\tau:=\ov{\omega}^{\,\tau}\backslash\,\{0\}$,
$\omega^\tau_M:=\{t_m\}_{m=1}^M$.
We also define the backward difference quotient, the mean value with the weight   $\sigma$ (independent of the meshes) and the backward shift in $t$
\[
 \ov{\pa}_t\Phi^m:=
 \frac{\Phi^m-\Phi^{m-1}}{\tau_m},\ \
 \Phi^{(\sigma)\,m}:=\sigma\Phi^m+(1-\sigma)\Phi^{m-1},\ \
 \check{\Phi}^m:=\Phi^{m-1}.
\]
\par We study the following finite-difference scheme, which is weighted in $t$ and averaged in $x$, on a finite mesh with an abstract approximate TBC for the initial-boundary value problem (\ref{sch})-(\ref{cc})
\begin{gather}
 C_\theta[\rho_h]\ov{\pa}_t U^m+{\mathcal A}_h U^{(\sigma)\,m}=F^m\ \
 \mbox{on}\ \ \omega_h\times\omega^\tau,
\label{dse}\\[1mm]
 U^m_0=g(t_m) \ \ \mbox{for}\ \ m\geq 1,
\label{ds0}\\[1mm]
 b_\infty\ov{\pa}_x U_J^{(\sigma)\,m}
 + h_J s_\theta^{-}(\rho_\infty\ov{\pa}_t U
 + c_\infty U^{(\sigma)})^m_J=b_\infty{\mathcal S}^m{\bf U}_J^m
 \ \ \mbox{for}\  m\geq 1,
\label{dsJ}\\[1mm]
 U^0=U^0_h\ \ \mbox{on}\ \ \ov{\omega}_h
\label{dsi}
\end{gather}
with the operator
${\mathcal A}_h W:=-\wh{\pa}_x\lf(b_h\ov{\pa}_x W \rg)+C_\theta [c_h] W$ and the functions
\[
 \rho_{hj}=\rho(x_{j-1/2}),\ b_{hj}=b(x_{j-1/2}),\ c_{hj}=c(x_{j-1/2}),\ F_j^m=f(x_j,t_m)
\]
and $U^0_{hj}=u^0(x_j)$ (for simplicity, for continuous $\rho$, $b$, $c$, $f$ and $u_0$). Thus $U^0_{hJ}=0$; we also assume that $U^0_{h0}=0$.
Here ${\mathcal S}^m$ is any linear operator acting in the space of functions given on the mesh $\omega_m^{\tau}\cup\{0\}$, and
${\bf U}_J^m:=\lf\{U_J^0,\dots,U_J^m\rg\}$.
\par Now we discuss the approximate TBC, i.e., the boundary condition (\ref{dsJ}). Let an equation
\begin{equation}
 \overset{\circ}{\pa}_x \lf[U^{(\sigma)}-\frac{\theta h^2}{b_\infty}\lf(\rho_\infty \ov{\pa}_t U+ c_\infty U^{(\sigma)}\rg)\rg]^m_J
 ={\mathcal S}^m {\bf U}_J^m\ \ \mbox{for}\ \ m\geq 1
\label{dtbc}
\end{equation}
serve as an (abstract) approximate TBC for (\ref{tbc}) at the node $x_J$; here we have discretized $\frac{\pa u}{\pa x}$ with weight in $t$ and symmetrically in $x$. We first write down equation (\ref{dse}) on the mesh $\omega_h\cup\{x_J\}$ and apply it at the node $x_J$ only in order to eliminate the values
$U^{(\sigma)}_{J+1}$ involved in the left-hand side of (\ref{dtbc}).
Namely, since
$\overset{\circ}{\pa}_xW_J=\ov{\pa}_xW_J+\frac{h}{2}\wh{\pa}_x\,\ov{\pa}_xW_J$,
taking into account (\ref{cc}) we get the boundary condition (\ref{dsJ}).
Importantly, the boundary condition (\ref{dsJ}) for $\mathcal{S}=0$ is the natural approximation of the Neumann boundary condition for this finite-difference scheme for $x=X$.
Such an approach was implemented earlier in \cite{DZ06,Z07,DZZ09}; it reliably leads to the computationally stable form of discrete TBCs (in contrast to some other approaches).
\par The corresponding three-point system of mesh equations for a vector $\{U_j^m\}_{j=0}^J$ of the solution values on the upper level, has the form:
\begin{gather*}
    \alpha_{\sigma,j}^mU_{j-1}^m+(\beta_{\sigma,j}^m+\beta_{\sigma,j+1}^m) U_j^m+\alpha_{\sigma,j+1}^mU_{j+1}^m
\\[1mm]
 =\alpha_{\sigma-1,j}^mU_{j-1}^{m-1}+(\beta_{\sigma-1,j}^m+\beta_{\sigma-1,j+1}^m) U_j^{m-1}
 +\alpha_{\sigma-1,j+1}^mU_{j+1}^{m-1}+h_{j+1/2}F_j^m
\end{gather*}
for $1\leq j\leq J-1$,
\begin{gather}
 U_0^m=g(t_m),\ \
 \alpha_{\sigma,J}^mU_{J-1}^m+(2\beta_{\sigma-1,J}^m-\mu_{m0})U_J^m
 \nn\\[1mm]
 =\alpha_{\sigma-1,J}^mU_{J-1}^{m-1}+(2\beta_{\sigma-1,J}^m-\mu_{(m-1)1})U_J^{m-1}
 +\sum\limits_{0\leq l\leq m-2} \mu_{l(m-l)}U_J^l,
\label{syst}
\end{gather}
compare with \cite{DZZ09}. Here the coefficients are given by formulas
\[
 \alpha_{\sigma,j}^m
 =\theta\lf(\frac{h_j\rho_{hj}}{\tau_m}+\sigma h_j c_{hj}\rg)
 -\sigma\frac{b_{hj}}{h_j},\
 \beta_{\sigma,j}^m
 =(\thalf-\theta)\lf(\frac{h_j\rho_{hj}}{\tau_m}
 +\sigma h_jc_{hj}\rg)
 +\sigma\frac{b_{hj}}{h_j}
\]
for $1\leq j\leq J$; in particular, for $j=J$, these expressions become more simple
\[
 \alpha_{\sigma,J}^m=\theta\lf(\frac{h\rho_{\infty} }{\tau_m}+\sigma hc_{\infty}\rg)
 -\sigma\frac{b_\infty}{h},\
 \beta_{\sigma,J}^m=(\thalf-\theta)\lf(\frac{h\rho_{\infty}}{\tau_m}
 +\sigma h c_{\infty}\rg)
 +\sigma\frac{b_\infty}{h}.
\]
Equation (\ref{syst}) is written assuming that the operator in the approximate TBC has the form ${\mathcal S}^m{\bf U}_J^m=\sum\limits_{l=1}^m \mu_{l(m-l)}U_J^l$.
\section{Stability of the finite-difference scheme on the finite mesh with the approximate TBC}
\label{sect3}
We consider the stability problem for the finite-difference scheme
(\ref{dse})-(\ref{dsi}) with respect to the initial data $U^0_h$, the free term $F$ and a perturbation in the boundary condition (\ref{dsJ}) and take $g(t)=0$. We need to introduce several mesh counterparts of the $L^2(\Omega)$-inner products
\begin{gather*}
 \lf(V,W\rg)_{\omega_h}=\sum_{j=1}^{J-1}V_jW_jh_{j+1/2},\ \
 \lf(V,W\rg)_{\wt{\omega}_h}=\sum_{j=1}^JV_jW_jh_j,
\\
 \lf(V,W\rg)_{\ov{\omega}_h}= \lf(V,\, W\rg)_{\omega_h} +V_JW_J\frac{h}{2}
\end{gather*}
and the corresponding norms
$\|\cdot\|_{\omega_h}$, $\|\cdot\|_{\wt{\omega}_h}$, $\|\cdot\|_{\ov{\omega}_h}$
(of course, for mesh functions given on $\omega_h$ or belonging to
 $H_0(\ov{\omega}_h)$).
\par We introduce a bilinear form
\[
 (U,W)_{C_{\theta} [\varkappa_h]}:=(C_{\theta} [\varkappa_h]U,W)_{\omega_h}
 +\varkappa_{hJ}(s_\theta^-U)_JW_J h_J\ \ \text{for}\ U,W\in H_0(\ov{\omega}_h).
\]
According to \cite{DZZ09}, it is symmetric for $\theta\leq\frac14$ and generates a norm $\|W\|_{C_{\theta}[\varkappa_h]}=(W,W)_{C_{\theta}[\varkappa_h]}^{1/2}$ for $\varkappa_h>0$ (or a seminorm for $\varkappa_h \geq 0$). Moreover, an inequality
\begin{equation}
 \sqrt{c_{\theta}\un{\rho}}\,\|W\|_{\overline{\omega}_h}
 \leq\|W\|_{C_{\theta}[\rho_h]}
 \leq\sqrt{(1+4\max\{-\theta,0\})\ov{\rho}}\,\|W\|_{\overline{\omega}_h}
\label{ineqn}
\end{equation}
holds for any $W\in H_0(\ov{\omega}_h)$ and $\theta\leq \frac14$,
where $c_{\theta}=1-4\max\{\theta,0\}$ and $\ov{\rho}=\max\limits_{1\leq j\leq J} \rho_{hj}$.
\par We also introduce mesh counterparts of the norms in $L^2(0,t_M)$ and $L^2(\Omega\times(0,t_M))$:
\[
 \|\Phi\|_{\omega^\tau_M}:=\Bigl(\sum_{m=1}^M(\Phi^m)^2\tau_m\Bigr)^{1/2},\
 \|\cdot\|_{\omega^{h,\tau}_M}:=\|\,\|\cdot\|_{\omega_h}\|_{\omega^\tau_M},\
 \|\cdot\|_{\wt{\omega}^{h,\tau}_M}:=\|\,\|\cdot\|_{\wt{\omega}_h}\|_{\omega^\tau_M}
\]
and also set
$\|\cdot\|_{C_{\theta}[\varkappa_h]^{\tau}_M}
 :=\|\,\|\cdot\|_{C{_\theta}[\varkappa_h]}\|_{\omega^\tau_M}$
for $\varkappa_h\geq 0$.
\begin{proposition}
\label{p1}
 Let $U$ be a solution to the finite-difference scheme (\ref{dse})-(\ref{dsi}) with a generalized boundary condition (\ref{dsJ}):
\begin{equation}
 b_\infty\ov{\pa}_x U_J^{(\sigma)\,m}+h_J s_\theta^{-}(\rho_h \ov{\pa}_t U
 + c_h U^{(\sigma)})^m_J=b_\infty{\mathcal S}^m{\bf U}_J^m+G^m
 \ \ \mbox{for}\ m\geq 1,
\label{dsJm}
\end{equation}
where $G$ is given on $\omega^\tau$.
Let the operator ${\mathcal S}$ satisfy an inequality
\begin{equation}
 \sum_{m=1}^M
 \lf({\mathcal S}^m\bPhi^m\rg)\Phi^{(\sigma)\,m}\tau_m
 \leq 0\ \ \mbox{for any}\ M\geq 1
\label{cs}
\end{equation}
for any function $\Phi$ given on $\ov{\omega}^{\,\tau}$ such that $\Phi^0=0$, where $\bPhi^m=\{\Phi^1,\ldots,\Phi^m\}$. Then, for $\sigma\geq\frac12$ and $\theta<\frac14$, the first energy bound
\begin{gather}
 \max\Bigl\{\,\max\limits_{0\leq m\leq M} \|U^m\|_{C_{\theta}[\rho_h]},\sqrt{2}\,\|U\|_{(1)}\Bigr\}
 \leq\|U_h^0\|_{C_{\theta}[\rho_h]}
\nn\\[1mm]
 +\frac{K_\sigma }{\sqrt{c_{\theta}\un{\rho}}}
 \sum\limits_{m=1}^M\|F^{(0)\,m}\|_{\omega_h}\tau_m
 +\sqrt{\frac{2}{\nu}}\,(\|F^{(1)}\|_{\wt{\omega}^{h,\tau}_M}
 +\sqrt{X}\|G\|_{\omega^\tau_M})
\label{sb}
\end{gather}
holds for any $M\geq 1$
and for any decomposition $F=F^{(0)}+\wh{\pa}_xF^{(1)}$ such that $F^{(1)}|_{j=J}=0$, with $K_\sigma:=2(\sigma+|1-\sigma|)$.
Here the norm $\|U\|_{(1)}$ is such that
\[
 \|U\|_{(1)}^2
 =(\sigma-\thalf)\|\sqrt{\tau}\,\ov{\pa}_tU\|_{C_{\theta}[\rho_h]^{\tau}_M}^2
 +\|\sqrt{b_h}\,\ov{\pa}_x U^{(\sigma)}\|_{\wt{\omega}^{h,\tau}_M}^2
 +\|U^{(\sigma)}\|_{C_{\theta}[c_h]^{\tau}_M}^2.
\]
\par The bound holds also in the case $\theta=\frac14$ provided that $F^{(0)}=0$ (one has to drop the summand with $F^{(0)}$).
\par Consequently, for $\sigma\geq\frac12$ and $\theta\leq \frac14$, the scheme has a unique solution.
\end{proposition}
\par {\bf Proof}. We take the $(\cdot\,,\cdot)_{\omega_h}$--inner product of equation (\ref{dse}) and a function $W\in H_0(\ov{\omega}_h)$, sum the result by parts (using the second assumption (\ref{cc})) and obtain
\begin{gather}
\lf(C_\theta[\rho_h]\ov{\pa}_t U^m,W\rg)_{\omega_h}
 +\lf(b_h\ov{\pa}_x U^{(\sigma)\,m},\ov{\pa}_xW\rg)_{\wt{\omega}_h}
 +\lf(C_\theta[c_h] U^{(\sigma)\,m}, W\rg)_{\omega_h}
\nonumber\\[1mm]
=b_\infty (\ov{\pa}_xU^{(\sigma)\,m}_J){W^m_J}
 +\lf(F^m, W\rg)_{\omega_h}.
\label{sumid}
\end{gather}
Choosing $W=U^{(\sigma)\,m}$ and applying the boundary condition
(\ref{dsJm}) and other  assumptions (\ref{cc}), we get
\begin{gather*}
(\ov{\pa}_t U^m, U^{(\sigma)\,m})_{C_{\theta}[\rho_h]}
 +(b_h\ov{\pa}_x U^{(\sigma)\,m},\ov{\pa}_xU^{(\sigma)\,m})_{\wt{\omega}_h}
 +(U^{(\sigma)\,m}, U^{(\sigma)\,m})_{C_{\theta}[c_h]}
\\[1mm]
 -b_\infty\lf({\mathcal S}^m{\bf U}_J^m\rg)U^{(\sigma)\,m}_J
 =(F^m,U^{(\sigma)\,m})_{\omega_h}+\,G^mU^{(\sigma)\,m}_J\ \
 \mbox{for}\ \ m\geq 1.
\end{gather*}
We multiply the result by $\tau_m$ and sum up it over $m=1,\dots,M$. Applying the formula
$U^{(\sigma)}=U^{(1/2)} + (\sigma - \thalf)\tau \ov{\pa}_tU$, we obtain {\it the first energy equality}
\begin{gather}
\ds{\thalf\|U^M\|^2_{C_{\theta}[\rho_h]}
 +(\sigma-\thalf)
 \|\sqrt{\tau}\,\ov{\pa}_tU\|_{C_{\theta}[\rho_h]^{\tau}_M}^2
 +\|\sqrt{b_h}\,\ov{\pa}_x U^{(\sigma)}\|_{\wt{\omega}^{h,\tau}_M}^2}
\nonumber\\[1mm]
 +\|U^{(\sigma)}\|_{C_{\theta}[c_h]^{\tau}_M}^2
 -b_{\infty}
     \sum\limits_{m=1}^M\lf({\mathcal S}^m{\bf U}_J^m\rg)U^{(\sigma)\,m}_J\tau_m
     =\thalf\|U_h^0\|^2_{C_{\theta}[\rho_h]}+I^{(1)M}
\label{p11}
\end{gather}
for $M\geq 1$, where
\[
 I^{(1)M}:=\sum\limits_{m=1}^M\lf[(F^m,U^{(\sigma)\,m})_{\omega_h}
 +G^mU^{(\sigma)\,m}_J\rg]\tau_m.
\]
\par For $F=F^{(0)}+\wh{\pa}_xF^{(1)}$, we sum the result by parts and derive the following bound
\begin{gather*}
 I^{(1)M}\leq\sum\limits_{m=1}^M\|F^{(0)\,m}\|_{\omega_h}
 \lf(|\sigma|\|U^m \|_{\omega_h}
 +|1-\sigma|\|U^{m-1} \|_{\omega_h}\rg)\tau_m
\\[1mm]
 +(\|F^{(1)}\|_{\wt{\omega}^{h,\tau}_M}+\sqrt{X}\|G\|_{\omega^\tau_M})
     \|\ov{\pa}_xU^{(\sigma)}\|_{\wt{\omega}^{h,\tau}_M}
\\[1mm]
 \leq (|\sigma|+|1-\sigma|)\sum\limits_{m=1}^M\|F^{(0)\,m}\|_{\omega_h}\tau_m
 \max\limits_{0\leq m\leq M}\|U^m\|_{\omega_h}
\\[1mm]
 +\frac{1}{\sqrt{\nu}}
 \lf(\|F^{(1)}\|_{\wt{\omega}^{h,\tau}_M}+\sqrt{X}\|G\|_{\omega^\tau_M}\rg)
 \|\sqrt{b_h}\,\ov{\pa}_xU^{(\sigma)}\|_{\wt{\omega}^{h,\tau}_M}.
\end{gather*}
Using the left-hand inequality (\ref{ineqn}), conditions $\sigma\geq\frac12$ and (\ref{cs}) and applying the standard argument
lead from the first energy equality (\ref{p11}) to bound (\ref{sb}). It is well-known that such a bound implies the existence and uniqueness of a solution to the finite-different scheme.
\begin{remark}
In the case $\theta=\frac14$, one can generalize bound (\ref{sb}) and next stability bounds for $F^{(0)}\not\equiv 0$ as well, see \cite{DZZ09}.
\end{remark}
\par We also define a symmetric bilinear form
\[
\mathcal{L}_{\omega_{h\theta}}(U,W)=(b_h\ov{\pa}_xU,\ov{\pa}_xW)_{\wt{\omega}_h}
+(C_\theta[c_h]U,W)_{\omega_h}
+c_\infty (s_\theta^{-}U_J)W_Jh_J
\]
for $U,W\in H_0(\ov{\omega}_h)$ and derive stability in the norm $\|W\|_{\mathcal{L}_{\omega_{h\theta}}}=\mathcal{L}_{\omega_{h\theta}}^{1/2}(W,W)$.
\begin{proposition}
\label{p1A}
Let $U$ be a solution to the finite-difference scheme
(\ref{dse})-(\ref{dsi}) with the generalized boundary condition (\ref{dsJm}) instead of (\ref{dsJ}).
Let the operator  ${\mathcal S}$ satisfy an inequality
\begin{equation}
 \sum_{m=1}^M
 \lf({\mathcal S}^m\bPhi^m\rg)\ov{\pa}_t\Phi^m\tau_m \leq 0\ \
 \mbox{for any}\ M\geq 1
\label{csA}
\end{equation}
for any function $\Phi$ given on $\ov{\omega}^{\,\tau}$ such that $\Phi^0=0$.
Then, for $\sigma\geq\frac12$ and $\theta<\frac14$, the second energy bound
\begin{gather}
\max\Bigl\{\,\max\limits_{0\leq m\leq M} \|U^m\|_{\mathcal{L}_{\omega_{h\theta}}},\sqrt{2}\,\|U\|_{(2)}\Bigr\}
\leq\|U_h^0\|_{\mathcal{L}_{\omega_{h\theta}}}
+\sqrt{\frac{2}{c_{\theta}\un{\rho}}}\,
\|F^{(0)}\|_{\omega^{h,\tau}_M}
\nonumber\\[1mm]
 +\frac{4}{\sqrt{\nu}}\Bigl[\|F^{(1)\,0}\|_{\wt{\omega}_h}
 +\sum\limits_{m=1}^M\|\ov{\pa}_tF^{(1)}\|_{\wt{\omega}_h}\tau_m
 +\sqrt{X}\Bigl(|G^0|+\sum\limits_{m=1}^M |\ov{\pa}_tG^m|\tau_m\Bigr)\Bigr]
\label{sbA}
\end{gather}
holds for any $M\geq 1$ and any decomposition
$F=F^{(0)}+\wh{\pa}_xF^{(1)}$ with $F^{(1)}|_{j=J}=0$.
Here
\[
 \|U\|_{(2)}^2
 =\sum\limits_{m=1}^M\Bigl[(\sigma-\thalf)\tau_m
     \|\ov{\pa}_tU^m\|_{\mathcal{L}_{\omega_{h\theta}}}^2     +\|\ov{\pa}_t U^m\|_{C_{\theta}[\rho_h]}^2\Bigr]\tau_m.
\]
\par The bound holds also in the case $\theta=\frac14$ provided that $F^{(0)}=0$.
\par Consequently, for $\sigma\geq\frac12$ and $\theta\leq \frac14$, the scheme has a unique solution.
\end{proposition}
\par {\bf Proof}.
We choose $W=\ov{\pa}_t U^m$ in (\ref{sumid}), apply the boundary condition
(\ref{dsJm}) and assumptions (\ref{cc}) and get
\begin{gather*}
\lf(\ov{\pa}_t U^m,\ov{\pa}_tU^m\rg)_{C_{\theta}[\rho_h]}
 + \mathcal{L}_{\omega_{h\theta}}(U^{(\sigma)\,m},\ov{\pa}_tU^m)
\\[1mm]
 -b_\infty\lf({\mathcal S}^m{\bf U}_J^m\rg)\ov{\pa}_t U^m_J
 =\lf(F^m, \ov{\pa}_t U^m_J \rg)_{\omega_h}+\,G^m \ov{\pa}_t U^m_J\ \
\mbox{for}\ \ m\geq 1.
\end{gather*}
\par We multiply the equality by $\tau_m$ and sum up it over $m=1,\dots,M$. Applying again the formula
$U^{(\sigma)}=U^{(1/2)} + (\sigma - \thalf)\tau \ov{\pa}_tU$, we obtain {\it the second energy equality}
\begin{gather}
 \sum\limits_{m=1}^M\lf\|\ov{\pa}_tU^m\rg\|^2_{C_\theta[\rho_h]}\tau_m
 +\thalf\|U^M\|^2_{\mathcal{L}_{\omega_{h\theta}}}
 +(\sigma-\thalf)
 \ds{\sum\limits_{m=1}^M
 \|\ov{\pa}_tU^m\|_{\mathcal{L}_{\omega_{h\theta}}}^2\tau_m^2}
\nonumber\\[1mm]
\ds{-b_\infty\sum_{m=1}^M
 \lf({\mathcal S}^m{\bf U}_J^m\rg)\ov{\pa}_t U^m_J\tau_m
 =\thalf\lf\|U_h^0\rg\|^2_{\mathcal{L}_{\omega_{h\theta}}}+I^{(2)M}}
\label{p11A}
\end{gather}
for $M\geq 1$, where
\[
 I^{(2)M}:=\sum\limits_{m=1}^M\lf[(F^m,\,\ov{\pa}_tU^m)_{\omega_h}
 +G^m\ov{\pa}_tU^m_J\rg]\tau_m.
\]
\par For $F=F^{(0)}+\wh{\pa}_xF^{(1)}$, we sum the result by parts in $t$ and $x$ and get
\begin{gather*}
 I^{(2)M}=\sum\limits_{m=1}^M(F^{(0)\,m},\,\ov{\pa}_tU^m)_{\omega_h}\tau_m
 -\lf.(F^{(1)\,m},\ov{\pa}_xU^m)_{\wt{\omega}_h}\rg|_{m=0}^{m=M}
\\[1mm]
 +\sum_{m=1}^M(\ov{\pa}_tF^{(1)\,m},\ov{\pa}_xU^{m-1})_{\wt{\omega}_h}\tau_m
 +\lf(G^mU_J^m)\rg|_{m=0}^{m=M}
 -\sum_{m=1}^M(\ov{\pa}_tG^m)U_J^{m-1}\tau_m
\\[1mm]
\leq\|F^{(0)}\|_{\omega^{h,\tau}_M}\|\ov{\pa}_tU\|_{\omega^{h,\tau}_M}
 +2\Bigl[\|F^{(1)\,0}\|_{\wt{\omega}_h}
\\[1mm]
 +\sum\limits_{m=1}^M\|\ov{\pa}_tF^{(1)}\|_{\wt{\omega}_h}\tau_m
 +\sqrt{X}\Bigl(|G^0|
 +\sum\limits_{m=1}^M |\ov{\pa}_t G^m|\tau_m\Bigr)\Bigr]
 \max\limits_{0\leq m\leq M}\|\ov{\pa}_xU^m\|_{\wt{\omega}_h}.
\end{gather*}
Using the left-hand inequality (\ref{ineqn}),
conditions $\sigma\geq\frac12$ and (\ref{csA}) and applying the standard argument
lead from the second energy equality (\ref{p11A}) to bound (\ref{sbA}).
\section{Stability of the finite-difference scheme on an infinite mesh}
\label{sect4}
In order to construct and study the discrete TBC, we first turn to the finite-difference scheme on an infinite mesh for the original problem (\ref{sch})-(\ref{ic}) on the half-axis
\begin{gather}
 C_\theta[\rho_h]\,\ov{\pa}_tU+{\mathcal A}_h U^{(\sigma)}=F
 \ \ \mbox{on}\ \
 \omega_{h,\infty}\times\omega^\tau,
\label{cne}\\[1mm]
 U^m_0=g(t_m) \ \ \mbox{for}\ m\geq 1,
\label{cnb}\\[1mm]
 U^0=U^0_h \ \ \mbox{on}\ \ \ov{\omega}_{h,\infty}.
\label{cni}
\end{gather}
Assumptions (\ref{cc}) are supposed to be fulfilled. Let $g(t)=0$ and $U^0_h|_{j=0}=0$.
\par We introduce the Hilbert spaces $H_h$ and $\wt{H}_h$ (mesh counterparts of $L^2(\mathbb{R}^+)$) consisting of functions $W$ given on the meshes respectively
$\ov{\omega}_{h,\infty}$ (and with $W_0=0$) and $\omega_{h,\infty}$
and such that $\|W\|_{\ell_2}^2=\sum\limits_{j=1}^\infty W_j^2<\infty$, equipped with the inner products
\[
 (V,W)_{H_h}:=\sum_{j=1}^\infty V_jW_jh_{j+1/2},\ \
 (V,W)_{\wt{H}_h}:=\sum_{j=1}^\infty V_jW_jh_j.
\]
Since $h_j=h$ for $j\geq J$, the conditions
$ \|W\|_{\ell_2}^2<\infty$,
$\|W\|_{H_h}^2<\infty$ and $\|W\|_{\wt{H}_h}^2<\infty$
are equivalent.
\par We define symmetric bilinear forms \cite{DZZ09}
\begin{gather*}
(V,W)_{C_{\theta}[\varkappa_h],\infty}:=(C_\theta[\varkappa_h]V,W)_{H_h},
\\[1mm]
 \mathcal{L}_ {H_{h\theta}}(V,W)=(b_h\ov{\pa}_xV,\ov{\pa}_xW)_{\wt{H}_h}
 +(C_\theta[c_h]V,W)_{H_h}
\end{gather*}
for $V,W\in H_h$. They generate norms $\|W\|_{C_{\theta}[\varkappa_h],\infty}=(W,W)_{C_{\theta}[\varkappa_h],\infty}^{1/2}$ for $\varkappa_h=\rho_h$ (a seminorm for $\varkappa_h=c_h$) and
$\|W\|_{\mathcal{L}_{H_{h\theta}}}=\mathcal{L}_{H_{h\theta}}^{1/2}(W,W)$. Moreover, an inequality
\[
 \sqrt{c_{\theta}\un{\rho}}\,\|W\|_{H_h}
 \leq\|W\|_{C_{\theta}[\rho_h],\infty}
 \leq\sqrt{(1+4\max\{-\theta,0\})\ov{\rho}}\,\|W\|_{H_h}
\]
holds for all $W\in H_h$ and $\theta\leq \frac14$ \cite{DZZ09}.
\par We also define the mesh counterparts of the norm in
$L^2(\mathbb{R}^+\times (0,t_M))$
\[
 \|\cdot\|_{H^{h,\tau}_M}:=\|\,\|\cdot\|_{H_h}\|_{\omega^\tau_M},\ \
 \|\cdot\|_{\wt{H}^{h,\tau}_M}:=\|\,\|\cdot\|_{\wt{H}_h}\|_{\omega^\tau_M}.
\]
\begin{proposition}
\label{p2}
Let $F=F^{(0)}+\wh{\pa}_xF^{(1)}$ with $F^{(0)\,m}\in H_h$ and $F^{(1)\,m}\in\wt{H}_h$
for any $m\geq 1$ and $U_h^0\in H_h$.
Then, for  $\sigma\geq 0$ and $\theta \leq \frac14$, there exists a unique solution $U^m\in H_h$, for all $m\geq 0$, to the finite-difference scheme (\ref{cne})-(\ref{cni}), and, for $\sigma\geq\frac12$ and $\theta<\frac14$, the first energy bound
\begin{gather}
\ds{\max\lf\{\,\max_{0\leq m\leq M}
\|U^m\|_{C_{\theta}[\rho_h],\infty},\sqrt{2}\,\|U\|_{(1),\infty}\rg\}}
\nn\\[1mm]
\ds{\leq\|U_h^0\|_{C_{\theta}[\rho_h],\infty}
+\frac{K_\sigma}{\sqrt{c_{\theta}\un{\rho}}}
\sum\limits_{m=1}^M\|F^{(0)\,m}\|_{H_h}\tau_m
 +\sqrt{\frac{2}{\nu}}\,\|F^{(1)}\|_{\wt{H}_M^{h,\tau}}}
\label{sb2}
\end{gather}
holds for any $M\geq 1$. Here
\[
 \|U\|_{(1),\infty}^2
 =\Bigl\|\sqrt{(\sigma-\thalf)\tau}\|\ov{\pa}_tU\|_{C_{\theta}[\rho_h] ,\infty}\Bigr\|_{\omega^\tau_M}^2
 +\|\|U^{(\sigma)}\|_{\mathcal{L}_{H_{h\theta}}}\|_{\omega^\tau_M}^2.
\]
\par The bound holds also in the case $\theta=\frac14$ provided that $F^{(0)}=0$.
\end{proposition}
\par{\bf Proof}. We extend ${\mathcal A}_h$ and $C_\theta[\rho_h]$ up to operators acting in $H_h$ by setting
$(\mathcal A_hW)_0:=0$ and $(C_\theta[\rho_h]W)_0:=0$.
By virtue of assumptions  (\ref{cc}) and the property $h_j=h$ for $j\geq J$, the operator $\mathcal A_h$ is bounded â $H_h$.
Moreover, $\mathcal A_h=\mathcal A_h^*>0$ since
\begin{equation}
 (\mathcal A_hW,V)_{H_h}=\mathcal{L}_ {H_{h\theta}}(V,W)
\label{sap}
\end{equation}
for any $W,V\in H_h$.
To establish equality (\ref{sap}), one can first transform the finite sum
$\sum\limits_{j= 1}^{j_1}({\mathcal A}_hW)_j V_jh_{j+1/2}$ by summing by parts (compare with the derivation of equality (\ref{sumid})) and then pass to the limit as $j_1\to\infty$
using the property $\lim_{j\to\infty}W_j=0$ for $W\in H_h$ (as in \cite{Z07} for $\theta=0$).
\par Now we rewrite equation (\ref{cne}), together with the homogeneous boundary condition  (\ref{cnb}), as an operator equation in $H_h$
\begin{equation}
 C_\theta[\rho_h]\,\ov{\pa}_t U+\mathcal A_h U^{(\sigma)}=F
 \ \ \mbox{on}\ \ \omega^\tau.
\label{cne1}
\end{equation}
Thus
\begin{equation}
 (C_\theta[\rho_h]+ \sigma\tau_m\mathcal A_h)U^m
 =(C_\theta[\rho_h]- (1-\sigma)\tau_m\mathcal A_h)\check{U}^m
 +\tau_m F^m
\label{cne2}
\end{equation}
for $m\geq 1$. For $\sigma\geq 0$, the operator $C_\theta[\rho_h]+ \sigma\tau_m\mathcal A_h$ is bounded, self-adjoint and positive definite and therefore invertible.
Since for $\check{U}^m\in H_h$ the right-hand side of equation (\ref{cne2}) also belongs to $H_h$, we find that the eqution has a unique solution $U^m \in H_h$.
\par Equation (\ref{cne1}) with the help of property (\ref{sap}) implies the first energy equality
\begin{gather}
\ds{\thalf\|U^M\|_{C_{\theta}[\rho_h],\infty}^2
 +\sum\limits_{m=1}^M\lf[(\sigma-\thalf)
 \lf\|\sqrt{\tau_m}\,\ov{\pa}_tU^m\rg\|_{C_{\theta}[\rho_h] ,\infty}^2
 +\|U^{(\sigma)\,m}\|_{C_{\theta}[c_h],\infty}^2\rg]\tau_m}
\nn\\[1mm]
 \ds{+\|\sqrt{b_h}\,\ov{\pa}_xU^{(\sigma)}\|_{\wt{H}^{h,\tau}_M}^2
 =\thalf\|U_h^0\|_{C_{\theta}[\rho_h],\infty}^2
 +\sum\limits_{m=1}^M(F^m, U^{(\sigma)\,m})_{H_h}\tau_m,}
\label{p11inf}
\end{gather}
compare with (\ref{p11}). Similarly to the proof of Proposition \ref{p1}, it implies bound (\ref{sb2}).
\begin{remark}
\label{r1p2}
Proposition \ref{p2} remains valid for any
$\rho_h$, $b_h$ and
$c_h$ satisfying the conditions $\rho_h\geq\un{\rho}>0$, $b_h\geq \nu>0$, $c_h\geq 0$ and
$\sup\limits_{j\geq 1}\lf(\rho_{hj}+b_{hj}+c_{hj}\rg)<\infty$.
\end{remark}
\begin{remark}
\label{r2}
The solution $U=U_\sigma$ to the finite-difference scheme (\ref{cne})-(\ref{cni}) (specified in Proposition \ref{p2}) depends continuously on $\sigma\geq\frac12$ for any $\theta\leq \frac14$.
\par Indeed, let $\sigma_0\geq 1/2$. Then the difference $U_{\sigma}-U_{\sigma_0}$
satisfies the operator equation in $H_h$
\begin{equation}
 C_\theta[\rho_h]\,\ov{\pa}_t (U_{\sigma}-U_{\sigma_0})
 +\mathcal A_h (U_{\sigma}-U_{\sigma_0})^{(\sigma)}
 =(\sigma_0-\sigma)\mathcal A_h(U_{\sigma_0}-\check{U}_{\sigma_0})\ \ \mbox{on}\ \ \omega^\tau.
\label{condep}
\end{equation}
Let $I_hV_j:=\sum_{1\leq k\leq j-1}V_kh_{k+1/2}$ for $j\geq 0$. Since $U_{\sigma}^0-U_{\sigma_0}^0=0$ and $A_hW=\wh{\pa}_x(-b_h\ov{\pa}_x W+I_hC_\theta [c_h] W)$ on $\omega_{h,\infty}$, the property follows from bound  (\ref{sb2}) (with $F^{(0)}=0$) applied to equation (\ref{condep}).
\end{remark}
\par We define the mesh counterparts of the norm in $L^2(X,\infty)$ such that
\begin{gather*}
 \|W\|_{D_h}^2:=\frac{h}{2}\,W_J^2+\sum_{j=J+1}^\infty W_j^2h,\ \
 \|W\|_{\wt{D}_h}^2:=\sum_{j=J+1}^\infty W_j^2h,
\\
 \|W\|_{s_{\theta},D_h}^2
 :=(s_\theta^+W)_JW_Jh+\sum_{j=J+1}^\infty (s_\theta W)_jW_jh;
\end{gather*}
concerning the correctness of the last definition, see  \cite{DZZ09}.
\begin{corollary}
\label{c2p2}
Let  $F_j^m=0$ and $U_{hj}^0=0$ for $j\geq J$ and $m\geq 1$. If the solution  $U^m\in H_h$, for all $m\geq 0$, to the scheme
(\ref{cne})-(\ref{cni}) satisfies the approximate TBC (\ref{dtbc}) with some operator $\mathcal{S}=\mathcal{S}_{\rm ref}$, then
an equality
\begin{gather*}
-b_\infty
 \sum\limits_{m=1}^M\lf(\mathcal{S}_{\rm ref}^m{\bf U}_J^m\rg)U^{(\sigma)\,m}_J\tau_m
 =\thalf \rho_{\infty}\lf\|U^M\rg\|_{s_{\theta},D_h}^2
\nonumber\\[1mm]
 +\sum\limits_{m=1}^M\lf[(\sigma-\thalf)\tau_m\rho_{\infty}
  \|\ov{\pa}_tU^m\|_{s_{\theta},D_h}^2
 +b_\infty\|\ov{\pa}_x U^{(\sigma)\,m}\|_{\wt{D}_h}^2
 +c_{\infty}\|U^{(\sigma)\,m}\|_{s_{\theta},D_h}^2\rg]\tau_m
\end{gather*}
holds for all $M\geq 1$. Its right-hand side is nonnegative for $\sigma\geq\frac12$.
\end{corollary}
\par {\bf Proof}. By virtue of equation
(\ref{cne}) at the node $x_J$ with $F_J^m=0$, relation
(\ref{dtbc}) is equivalent to the boundary condition (\ref{dsJ}); thus, the solution to the scheme (\ref{cne})-(\ref{cni}) satisfies the scheme (\ref{dse})-(\ref{dsi}) as well.
Taking the difference of the energy equalities (\ref{p11inf}) and (\ref{p11}) (with $G=0$) and applying simple identities
\begin{gather*}
 \|W\|_{H_h}^2=\|W\|_{\ov{\omega}_h}^2
 +\|W\|_{D_h}^2,\ \
 \|W\|_{\wt{H}_h}^2=\|W\|_{\wt{\omega}_h}^2+\|W\|_{\wt{D}_h}^2,
\\[1mm]
\|W\|_{C_{\theta}[\varkappa_h],\infty}^2=\|W\|_{C_{\theta}[\varkappa_h]}^2
 +\varkappa_{hJ}\|W\|_{s_{\theta},D_h}^2
\end{gather*}
for any $W\in H_h$ and $\vk_h=\rho_h,c_h$, we obtain the announced equality.
\par We also derive stability in another norm.
\begin{proposition}
\label{p2a}
Let $F=F^{(0)}+\wh{\pa}_xF^{(1)}$ with $F^{(0)\,m}\in H_h$ and $F^{(1)\,m}\in\wt{H}_h$ for any $m\geq 1$ and $U_h^0\in H_h$.
Then, for $\sigma\geq\frac12$ and $\theta<\frac14$, the second energy bound
\begin{gather}
 \max\Bigl\{\,\max\limits_{0\leq m\leq M} \|U^m\|_{\mathcal{L}_ {H_{h\theta}}},  \sqrt{2}\,\|U\|_{(2),\infty}\Bigr\}\leq\|U_h^0\|_{\mathcal{L}_{H_{h\theta}}}
\nonumber\\[1mm]
 +\sqrt{\frac{2}{c_{\theta}\un{\rho}}}
 \|F^{(0)\,m}\|_{H^{h,\tau}_M}
 +\frac{4}{\sqrt{\nu}}\Bigl(\|F^{(1)\,0}\|_{\wt{H}_h}
 +\sum\limits_{m=1}^M\|\ov{\pa}_tF^{(1)\,m}\|_{\wt{H}_h}\tau_m\Bigr)
\label{sbAA}
\end{gather}
holds for the solution  $U^m\in H_h$, for all $m\geq 0$, to the finite-difference scheme (\ref{cne})-(\ref{cni}) and any $M\geq 1$. Here
 \[
 \|U\|_{(2),\infty}^2
 =\Bigl\|\sqrt{(\sigma-\thalf)\tau}\|\ov{\pa}_tU\|_{\mathcal{L}_ {H_{h\theta}}}\Bigr\|_{\omega^\tau_M}^2
 +\|\|\ov{\pa}_t U\|_{C_{\theta}[\rho_h],\infty}\|_{\omega^\tau_M}^2.
\]
\par The bound holds also in the case $\theta=\frac14$ provided that $F^{(0)}=0$.
\end{proposition}
\par {\bf Proof}.
The second energy equality
\begin{gather}
 \sum\limits_{m=1}^M\|\ov{\pa}_tU^m\|_{C_{\theta}[\rho_h], \infty}^2\tau_m
 +\thalf\lf\|U^M\rg\|_{\mathcal{L}_{H_{h\theta}}}^2
 +\sum\limits_{m=1}^M(\sigma-\thalf)
     \|\ov{\pa}_tU^m\|_{\mathcal{L}_{H_{h\theta}}}^2\tau_m^2
\nn\\[1mm]
 =\thalf\|U_h^0\|_{\mathcal{L}_{H_{h\theta}}}^2
 +\sum\limits_{m=1}^M\lf(F^m,\,\ov{\pa}_tU^m\rg)_{H_h}\tau_m
\label{p11AA}
\end{gather}
holds, compare with (\ref{p11A}). Similarly to the proof of Proposition \ref{p1A}, it implies bound (\ref{sbAA}).
\begin{corollary}
\label{r1p21}
Under the hypotheses of Corollary \ref{c2p2}, an equality
\begin{gather*}
 -b_\infty
 \sum\limits_{m=1}^M\lf({\mathcal S}_{\rm ref}^m{\bf U}_J^m\rg)\ov{\pa}_tU^m_J\tau_m
 =\thalf\lf(b_\infty\|\ov{\pa}_xU^M\|_{\wt{D}_h}^2
 +c_\infty\|U^M\|_{s_\theta,D_h}^2\rg)
\\[1mm]
 +\sum\limits_{m=1}^M\lf[\rho_\infty\|\ov{\pa}_tU^m\|_{s_{\theta},D_h}^2
 +(\sigma-\thalf)\tau_m
 \lf(b_\infty\|\ov{\pa}_x\ov{\pa}_tU^m\|_{\wt{D}_h}^2
 +c_\infty\|\ov{\pa}_tU^m\|_{s_\theta,D_h}^2\rg)\rg]\tau_m
\end{gather*}
holds for any $M\geq 1$. Its right-hand side is nonnegative for $\sigma\geq\frac12$.
\end{corollary}
\par {\bf Proof}. The result is derived by taking the difference of (\ref{p11AA}) and (\ref{p11A}) (with $G=0$).
\par By definition, \textit{the discrete TBC} is an approximate TBC (\ref{dtbc}) with the operator $\mathcal{S}=\mathcal{S}_{\rm ref}$. It will be explicitly constructed in the next section.
Corollaries \ref{c2p2} and \ref{r1p21} clarify the energy meaning of conditions (\ref{cs}) and (\ref{csA}) for the discrete TBC, for $\Phi=U_J$, and are exploited below to prove the conditions (for any $\Phi$).
\section{Derivation and analysis of the discrete TBC}
\label{sect5}
Now we turn to derivation of the explicit form for the discrete TBC in the form
(\ref{dtbc}) and verification of inequalities (\ref{cs}) and (\ref{csA}) for it.
We confine ourselves by the case of the uniform mesh
$\ov{\omega}^{\,\tau}$, i.e., $\tau_m=\tau$ for $m\geq 1$.
Consider an auxiliary finite-difference problem on the uniform part of the infinite mesh in $x$
\begin{gather}
 \rho_\infty s_\theta\ov{\pa}_t U+{\mathcal A}_{h,\infty} U^{(\sigma)}=0\ \
 \mbox{on}\ \
 (\omega_{h,\infty}\backslash\,\omega_h)\times\omega^\tau,
\label{ae}\\[1mm]
 U|_{j=J}=\Phi,\ \ \mbox{with}\ \
 |\Phi|_{\infty,\,q_0}:=\sup\limits_{m\geq 0}q_0^{-m}|\Phi^m|<\infty,\ \ \Phi^0=0,
\label{ab}\\[1mm]
 U^0_j=0\ \ \mbox{for}\ \ j\geq J-1
\label{ai}
\end{gather}
for some $q_0\geq 1$. Here the limiting finite-difference operator
\[
 {\mathcal A}_{h,\infty}W
 :=-b_\infty\wh{\pa}_x\,\ov{\pa}_x W+c_{\infty}s_\theta W
 \ \ \mbox{on}\ \ \omega_{h,\infty}\backslash\,\omega_h.
\]
has appeared. We seek for the solution satisfying the following property
\begin{equation}
 \|U\|_{2, \infty,\, q}
 :=\sup\limits_{m\geq0}q^{-m}\Bigl(\sum_{j=J-1}^{\infty}|U_j^m|^2\Bigr)^{1/2}<\infty
\label{ttt}
\end{equation}
for sufficently large $q>q_0$.
\par Since the coefficients are constant and the meshes are uniform,
the stated problem can be solved explicitly. For a mesh function $\Phi$: $\overline{\omega}^{\,\tau} \to {\mathbb C}$ such that $|\Phi|_{\infty,\,q}<\infty$ for some $q>0$, recall \textit{the reproducing function}
\[
 \widetilde{\Phi}(z)\equiv {\mathcal T}[\Phi](z):=\sum_{m=0}^{\infty}\Phi^m z^m
 \in A(D_{1/q}),
\]
i.e., analytic in the disc
$D_{1/q}:=\{|z|<\frac{1}{q}\}\subset\mathbb{C}$, satisfying a bound
\begin{equation}
 |\widetilde{\Phi}(z)|\leq \frac{|\Phi|_{\infty,\,q}}{1-q|z|}\ \ \mbox{for}\ \ |z|<q^{-1}.
\label{bzt}
\end{equation}
Conversely, for a function $p\in A(D_r)$ (for some $r>0$), the transformation $\Phi={\mathcal T}^{-1}[p]$ such that
\begin{equation}
\Phi^m=\frac{p^{(m)}(0)}{m!}
=\frac{1}{2\pi} \int_0^{2\pi}
\left. \frac{p(z)}{z^m}\right|_{z=r_1e^{i\varphi}} d\varphi\ \ \mbox{for any}\ m\geq 0,\ 0<r_1<r,
\label{12}
\end{equation}
is well defined implying the Cauchy inequality
\begin{equation}
 |\Phi|_{\infty,\,1/r_1}\leq\max_{|z|=r_1}|p(z)|.
\label{cauchy}
\end{equation}
Hereafter $i$ is the imaginary unit, and $\Rea z$ and $\Ima z$ are real and imaginary parts of $z\in {\mathbb C}$.
\par Taking into account conditions (\ref{ai}) and (\ref{ttt}),
for $|z|<r=\frac1q$, we calculate
\begin{gather}
 \mathcal{T}\lf[\rho_\infty s_\theta\ov{\pa}_t U_j
 +{\mathcal A}_{h,\infty} U^{(\sigma)}_j\rg](z)
 =\rho_\infty s_\theta\,\frac{1-z}{\tau}\,\widetilde{U}_j(z)
+z^{(\sigma)}{\mathcal A}_{h,\infty}\widetilde{U}_j(z)
\nn\\[1mm]
 =\Bigl\{\rho_\infty\frac{\theta(1-z)}{\tau}
 +\Bigl(-\frac{b_\infty}{h^2}+\theta c_\infty\Bigr)z^{(\sigma)}\Bigr\}
 \bigl(\widetilde{U}_{j+1}(z)-2\gamma(z)\widetilde{U}_j(z)
 +\widetilde{U}_{j-1}(z)\bigr)
\nonumber\\[1mm]
 =-\frac{b_\infty}{h^2}\,d(z)\bigl(\widetilde{U}_{j+1}(z)-2\gamma(z)\widetilde{U}_j(z)
 +\widetilde{U}_{j-1}(z)\bigr)
\label{schz}
\end{gather}
provided that $d(z)\neq 0$, with $z^{(\sigma)}:=\sigma+(1-\sigma)z$.
Hereafter, for $j\geq J-1$, we extend $\lf.U_j^m\rg|_{m=-1}:=0$ so that
$\ov{\pa}_t U_j^0=U^{(\sigma)\,0}_j=0$.
The coefficients $\gamma(z)$ and $d(z)$ are expressed by formulas
\begin{gather*}
 d(z)=2a_1\theta(z-1)+(1-2a_0\theta)z^{(\sigma)},\ \ \gamma(z):=1+\frac{a_1(1-z)+a_0z^{(\sigma)}}{d(z)},
\\[1mm]
 a_1:=\frac{h^2\rho_\infty}{2\tau b_\infty}>0,\ \
 a_0:=\frac{h^2c_\infty}{2b_\infty}\geq 0.
\end{gather*}
\par By virtue of (\ref{ae}) and (\ref{schz}) a difference equation
\begin{equation}
 \widetilde{U}_{j+1}(z)-2\gamma(z)\widetilde{U}_j(z)
 +\widetilde{U}_{j-1}(z)=0\ \ \mbox{for}\ \ j\geq J
\label{rrp}
\end{equation}
holds. The corresponding characteristic equation has the form
\begin{equation}
 \nu^2(z)-2\gamma(z)\nu(z)+1=0.
\label{cheq}
\end{equation}
\par Notice that $d(0)=0$ for $\sigma=\sigma_0:=\frac{2a_1\theta}{1-2a_0\theta}$.
\begin{lemma}
\label{l1}
For $\sigma>0$, $\sigma\neq\sigma_0$ and $\theta\leq \frac14$, the quadratic equation (\ref{cheq}) has roots $\nu_1,\nu_2\in A(D_r)$, for sufficiently small $r>0$, such that
\[
 \nu_1(z)=Z_1^{-1}(\gamma(z)),\ \ 0<|\nu_1(z)|<1,\ \
 \nu_2(z)=Z_2^{-1}(\gamma(z))=\frac{1}{\nu_1(z)},\ \ |\nu_2(z)|>1,
\]
where $Z_k^{-1}(\gamma)=\gamma+(-1)^{k+1}\sqrt[*]{\gamma^2-1}$, $k=1,2$ are analytic branches of the two-valued inverse function to the elementary Zhukovskii function $Z(z)=\thalf(z+z^{-1})$ defined in $\mathbb{C}$ with the cross-cut along the segment $[-1,1]$ of the real axis \cite{LSh87}.
\end{lemma}
\par {\bf Proof}. The presented formulas are rather elementary.
The property $\nu_1,\nu_2\in A(D_r)$ holds provided that
$\gamma(D_r)\subset\mathbb{C}\setminus [-1,1]$.
For validity of the latter property for sufficiently small $r>0$, it is required that $|\gamma(0)|<\infty$ (i.e., $\sigma\neq\sigma_0$) and $|\gamma(0)|>1$.
For $\sigma\neq\sigma_0$ and $2a_0\theta\neq 1$, formulas
\begin{gather*}
 \gamma(0)-1=\frac{a_1+a_0\sigma}{(1-2a_0\theta)(\sigma-\sigma_0)},\
 \gamma(0)+1=\frac{(1-4\theta)a_1+[2+(1-4\theta)a_0]\sigma}
 {(1-2a_0\theta)(\sigma-\sigma_0)}
\end{gather*}
hold. Since $\sigma>0$ and $\theta\leq \frac14$, the nominators of the both formulas are positive and thus
$\gamma(0)>1$ for $d(0)=(\sigma-\sigma_0)(1-2a_0\theta)>0$, or
$\gamma(0)<-1$ for $d(0)=(\sigma-\sigma_0)(1-2a_0\theta)<0$.
If $2a_0\theta=1$, then once again $\gamma(0)=1-a_0-a_0^2\sigma/a_1<-1$.
\par The following result corresponds to Proposition 5.3 in \cite{DZZ09}.
\begin{proposition}
\label{pinf}
For $\sigma>0$, $\sigma\neq\sigma_0$ and $\theta\leq \frac14$, the solution to the problem (\ref{ae})-(\ref{ttt}) exists, is unique and is given by a formula
\begin{equation}
 U_j=\mathcal{T}^{-1}\lf[\nu_1^{j-J}(z)\wt{\Phi}(z)\rg]\ \
 \mbox{for}\ \ j\geq J-1.
\label{gfsol}
\end{equation}
\par This solution satisfies a bound
\begin{equation}
 \|U\|_{2,\infty,\,q}\leq C|\Phi|_{\infty,\,q_0}
\label{poslusl}
\end{equation}
for sufficiently large $q>q_0$. For real $\Phi$, it is real too.
\end{proposition}
\par {\bf Proof}.
Let $\sigma\neq\sigma_0$. By taking into account Lemma \ref{l1}, for $z\in D_r$ the general solution to the difference equation (\ref{rrp}) has the form
\[
 \widetilde{U}_j(z)=c_1(z)\nu_1^{j-J+1}(z)+c_2(z)\nu_2^{j-J+1}(z)\ \
 \mbox{for}\ \ j\geq J-1
\]
with any $c_1(z)$ and $c_2(z)$. By virtue of bounds (\ref{ttt}) and (\ref{bzt})
we find that $c_2(z)\equiv 0$, and then from condition (\ref{ab})
we derive a formula (taking into account that $\nu_1(z)\neq 0$)
\begin{equation}
 \widetilde{U}_j(z)=\widetilde{\Phi}(z)\nu_1^{j-J}(z)\ \
 \mbox{for}\ \ j\geq J-1.
\label{psiz}
\end{equation}
Since $\widetilde{\Phi}\nu_1^{j-J}\in A(D_r)$ exploiting Lemma \ref{l1}, if the solution to the problem (\ref{ae})-(\ref{ttt}) exists, then it is given by formula (\ref{gfsol}).
\par Conversely, the function given by formula (\ref{gfsol}) satisfies an equation
\[
 \mathcal{T}\lf[\rho_\infty s_\theta\ov{\pa}_t U_j
 +{\mathcal A}_{h,\infty} U^{(\sigma)}_j\rg](z)=0\ \ \mbox{for}\ \ j\geq J,\ z\in D_r
\]
and therefore equation (\ref{ae}) as well. It also satisfies conditions
 (\ref{ab}) and (\ref{ai}).
\par By virtue of (\ref{cauchy})
and (\ref{bzt}) we get that, for any $j\geq J$ and $m\geq 1$, bounds
\[
 q_1^{-m}|U_j^m|\leq \max\limits_{|z|=r_1}|\nu_1^{j-J}(z)\widetilde{\Phi}(z)|
 \leq C_{0j}^{j-J}\frac {|\Phi|_{\infty,\,q_0}}{1-q_0r_1}
\]
hold for $q_1=\frac{1}{r_1}>\frac{1}{r}>q_0$, with
\[
C_{0(J-1)}=\nu_{1\min}=\min\limits_{|z|=r_1}|\nu_1(z)|,\ \ C_{0j}=\nu_{1\max}=\max\limits_{|z|=r_1}|\nu_1(z)|\ \ \text{for}\ j\geq J.
\]
Therefore bound (\ref{poslusl}) holds with $C=\frac{C_1}{1-q_0r_1}$, where
\[
 C_1^2=\sum\limits_{j=J-1}^{\infty}C_{0j}^{2(j-J)}
 =\frac{1}{\nu_{1\min}^{2}} + \frac{1}{1-\nu_{1\max}^{2}},\ \ C_1>0.
\]
\par For real $\Phi$, the functions $\widetilde{\Phi}(z)$ and $\gamma(z)$ are real as well for $z\in\mathbb{R}$.
If in addition $|z|\leq r$, then $|\gamma(z)|>1$ and thus
$\nu_1(z)=Z_1^{-1}(\gamma(z))$ is real. Therefore $U$ is real too, see (\ref{12}).
\par The proofs of Lemma \ref{l1} and Proposition \ref{pinf} remain valid also for $\sigma=0$, $\theta<\frac14$ and $\theta\neq 0$.
\par We go back to the derivation of the discrete TBC. By virtue of formula (\ref{psiz}) we have
\[
 \widetilde{U}_{j+1}(z)-\widetilde{U}_{j-1}(z)
 =\lf(\nu_1(z)-\frac{1}{\nu_1(z)}\rg)\widetilde{U}_j(z)\ \ \mbox{for}\ \ j\geq J.
\]
Therefore it is easy to check that a formula
\begin{equation}
 \mathcal{T}\lf[\overset{\circ}{\pa}_x
 \Bigl(U^{(\sigma)}
 -\frac{\theta h^2}{b_\infty}\,(\rho_\infty \ov{\pa}_t U
 + c_\infty U^{(\sigma)})\Bigr)_J\rg]
 =\frac{1}{2h}\,d(z)(\nu_1-\nu_2)(z)\,\wt{\Phi}(z)
\label{zdlt}
\end{equation}
holds, compare with (\ref{schz}).
By virtue of the well known formula for the multipklication of two poer series
it leads to the discrete TBC (\ref{dtbc}) with the operator ${\mathcal S}={\mathcal S}_{\rm ref}$ of the discrete convolution form
\begin{equation}
 {\mathcal S}_{\rm ref}^m {\bPhi}^m
 :=\frac{1}{2h}\lf(R*\Phi\rg)^m
 =\frac{1}{2h}\,\sum_{q=0}^{m} R^q\Phi^{m-q}
 \ \ \mbox{for}\ \ m\geq 1,
\label{mdtbc}
\end{equation}
with the kernel
\begin{equation}
R:=\mathcal{T}^{-1}[d(z)(\nu_1-\nu_2)(z)].
\label{mdtbcr}
\end{equation}
\par Let us see that Propositions \ref{p1} and \ref{p1A} on stability are valid for ${\mathcal S}={\mathcal S}_{\rm ref}$.
\begin{proposition}
\label{p3}
For $\sigma\geq\frac12$, $\sigma\neq\sigma_0$ and $\theta\leq \frac14$, the operator ${\mathcal S}={\mathcal S}_{\rm ref}$ of the discrete TBC (\ref{mdtbc}) satisfies inequalities (\ref{cs}) and (\ref{csA}).
\end{proposition}
\par {\bf Proof}. We apply a method first suggested in \cite{DZ06}. Fix any
$M\geq 1$ and real values $\Phi^1,\dots,\Phi^M$.
Extend $\Phi^m=0$ for $m=-1,0$ and $m>M$.
We define a function $U$ by formula (\ref{gfsol}) for $j\geq J-1$
and set, for example, $U_j^m:=0$ for $0\leq j<J-1$ and $m\geq 0$,
and then set $F:=C_\theta[\rho_h]\ov{\pa}_t U+{\mathcal A}_h U^{(\sigma)}$.
\par By virtue of Proposition \ref{pinf}, the constructed function $U$ serves as the real solution to the problem (\ref{ae})-(\ref{ai}) and therefore as one to the scheme (\ref{cne})-(\ref{cni}), where
$F=0$ on $(\omega_{h,\infty}\backslash\,\omega_h)\times\omega^\tau$ and $U_h^0=0$.
Then Corollary \ref{c2p2} implies inequality (\ref{cs}) whereas Corollary \ref{r1p21} implies inequality (\ref{csA}) since $U_J^m=\Phi^m$ for $0\leq m\leq M$.
%
\par The case $\sigma=\sigma_0\geq\frac12$ will be also covered in Remark \ref{rem4} below.
\par We find the kernel $R$ explicitly. We introduce quantities by recurrence formulas
\begin{gather}
p_{m,\alpha,\beta}=\frac{2m-1}{m}\,\b p_{m-1,\alpha,\beta}
 -\frac{m-1}{m}\,\a p_{m-2,\alpha,\beta}\ \ \mbox{for}\ m\geq 1,
\label{pmvk}\\[1mm]
p_{0,\alpha,\beta}=1,\ \ p_{m,\alpha,\beta}=0\ \ \mbox{for}\ \ m<0
\label{pmvk0}
\end{gather}
with parameters $\alpha$, $\beta$. Their close connecton to the classical Legendre polynomials will be shown below.
\begin{proposition}
\label{p4}
For $\sigma>0$, $\sigma\neq\sigma_0$ or $\sigma=\sigma_0\geq\frac12$, and $\theta\leq \frac14$, the kernel $R$ of the operator of the discrete TBC (\ref{mdtbc}) is real and is expressed by an explicit formula
\begin{equation}
 R^m=2a_1\sqrt{\de}\,
 \frac{1}{2m-1}\,\lf[p_{m,\alpha,\beta}-\a p_{m-2,\alpha,\beta}\rg]\ \
 \mbox{for}\ \ m\geq 0,
\label{formr}
\end{equation}
with the quantities $\a$, $\b$ and $\de$ of the form
\begin{gather}
\a=\a_0\a_1,\ \ \beta=\frac{\a_0+\a_1}{2},
\label{a0a1n}\\[1mm]
\a_0=1-\frac{d_0}{1+\sigma d_0},\ \
\a_1=1-\frac{d_0(1-4\theta)+d_1}
          {(1+\sigma d_0)(1-4\theta)+\sigma d_1},
\label{a0a1n1}\\[1mm]
 \de:=(1+\sigma d_0)[(1+\sigma d_0)(1-4\theta)+\sigma d_1]>0,
\label{fde}\\[1mm]
 d_0:=\frac{a_0}{a_1}=\frac{c_\infty}{\rho_\infty}\,\tau\geq 0,\ \
 d_1:=\frac{2}{a_1}=4\frac{b_\infty}{\rho_\infty}\,\frac{\tau}{h^2}>0.
\label{fded}
\end{gather}
\end{proposition}
\par {\bf Proof}. Let first $\sigma\neq\sigma_0$.
By Lemma \ref{l1} for $z\in D_r$ we have
\begin{equation}
 (\nu_1-\nu_2)(z)=2\sqrt[*]{\gamma^2(z)-1}.
\label{n1n2}
\end{equation}
One can straightforwardly verify that
\begin{gather}
 d^2(z)\lf(\gamma^2(z)-1\rg)=
\nn\\[1mm]
 =a_1^2(1-z+d_0z^{(\sigma)})[(1-4\theta)(1-z)+(d_1+d_0(1-4\theta))z^{(\sigma)}]=
\nn\\[1mm]
 =a_1^2\delta(\alpha_0z-1)(\alpha_1z-1)
 =a_1^2\delta(\alpha z^2-2\beta z+1),
\label{g2zm1}
\end{gather}
where the coefficients $d_0$, $d_1$, $\a_0$, $\a_1$ and $\a$, $\b$, $\de$ are given by formulas (\ref{a0a1n})-(\ref{fded}).
\par Let first also $\a\neq 0$. We write down a formula
\[
 \alpha z^2-2\beta z+1
 =(\varkappa z)^2-2\mu\varkappa z+1,
\]
where $\vk=\sqrt{\a}$ and $\mu=\frac{\b}{\sqrt{\a}}$ are real for $\a>0$, or
$\vk=i\sqrt{|\a|}$ and $\mu=-i\frac{\b}{\sqrt{|\a|}}$ are purely imaginary for $\a<0$;
herewith $\vk^2=\a$ and $\vk\mu=\b$ are always real.
Inserting the formula into (\ref{g2zm1}) leads to
\begin{equation}
 d(z)\sqrt[*]{\gamma^2(z)-1}
 =-a_1\sqrt{\de}
 \sqrt[+]{(\varkappa z)^2-2\mu\varkappa z+1}
\label{g2z}
\end{equation}
at least for sufficiently small $z$ (in accordance with the proof of Lemma \ref{l1}), where $\sqrt[+]{w}$ is an analytic branch of
$\sqrt{w}$ on $\mathbb{C}$ with the cross-cut along the negative real half-axis $\Rea w<0$ such that $\sqrt[+]{1}=1$.
Formulas (\ref{n1n2}) and (\ref{g2z}) imply an equality
\[
 d(z)(\nu_1-\nu_2)(z)=-2a_1\sqrt{\de}\,
 \frac{(\varkappa z)^2-2\mu\varkappa z+1}
      {\sqrt[+]{(\varkappa z)^2-2\mu\varkappa z+1}}.
\]
The following generalized formula for the reproducing function of the Legandre polynomials holds \cite{DZ06}:
\[
 \sum_{m=l}^\infty\varkappa^mP_{m-l}(\mu)z^m
 =\frac{(\varkappa z)^{l}}{\sqrt[+]{(\varkappa z)^2-2\mu\varkappa z+1}}
\]
for any $\varkappa\in \mathbb{C}$, integer $l\geq 0$ and
sufficiently small $z$, that easily follows from the classical one in the case
$\vk=1$, $l=0$ \cite{LSh87}. Therefore
\[
 d(z)(\nu_1-\nu_2)(z)
 =-2a_1\sqrt{\de}
 \sum_{m=0}^\infty\varkappa^m\lf[P_{m-2}(\mu)-2\mu P_{m-1}(\mu)+P_m(\mu)\rg]z^m,
\]
where $P_m(\mu)\equiv 0$ for $m<0$. Applying the recurrence relation for the Legendre polynomials
\begin{equation}
 \mu P_{m-1}(\mu)=\frac{m-1}{2m-1}\,P_{m-2}(\mu)+\frac{m}{2m-1}\,P_m(\mu)\ \
 \mbox{for}\ \ m\geq 0,
\label{legrec}
\end{equation}
one can simplify the last formula as follows
\[
 d(z)(\nu_1-\nu_2)(z)
 =2a_1\sqrt{\de}
 \sum_{m=0}^\infty\frac{\varkappa^m}{2m-1}\,(P_m-P_{m-2})(\mu)z^m,
\]
i.e., to derive a formula
\begin{equation}
 R^m=2a_1\sqrt{\de}\frac{\varkappa^m}{2m-1}\,(P_m-P_{m-2})(\mu),\ \ m\geq 0.
\lb{formr1}
\end{equation}
\par Following \cite{ME05,Z07}, we introduce modified Legendre polynomials
$p_{m,\vk}(z):=\vk^mP_m(z)$.
From (\ref{legrec}) clearly $p_{m,\alpha,\beta}=p_{m,\vk}(\mu)$
satisfy recurrence equalities (\ref{pmvk}) and (\ref{pmvk0}) and,
in particular, they are real. Therefore formula (\ref{formr}) is proved.
\par For $\a=0$ formula (\ref{g2z}) is simplified and takes the form
$d(z)\sqrt[*]{\gamma^2(z)-1} =-a_1\sqrt{\de} \sqrt[+]{1-2\beta z}$.
Hence one can easily verify that formulas
(\ref{pmvk}), (\ref{pmvk0}) and (\ref{formr}) remain valid and even are simplified in the case $\a=0$.
\par Owing to continuous dependence of $U$ on $\sigma\geq\frac12$ (see Remark \ref{r2}) and $R$ on $\sigma\geq\frac12$ (it is clear), one can pass to the limit as $\sigma\to\sigma_0$ on the left-hand side of the discrete TBC (\ref{dtbc}) with
${\mathcal S}={\mathcal S}_{\rm ref}$ and on the right-hand side of equality (\ref{mdtbc}). This justifies the validity of formula
(\ref{mdtbc}), for $R$ of the form (\ref{formr}), in the case $\sigma=\sigma_0\geq \frac12$ (that is possible only provided that $\frac{1}{4a_1+2a_0}\leq\theta<\frac{1}{2a_0}$).
\begin{remark}
\label{rem4}
Proposition \ref{p3} remains valid also in the case $\sigma=\sigma_0\geq\frac12$.
To see this, it suffices to insert formula (\ref{mdtbc}) into inequalities (\ref{cs}) and (\ref{csA}) for
${\mathcal S}={\mathcal S}_{\rm ref}$ that have been already proved, for $\sigma\neq\sigma_0$, and
pass to the limit as $\sigma\to\sigma_0$ taking into account the continuous dependence of $R$ on $\sigma$.
\end{remark}
\par In practical computations, recurrence equalities for $R$ are more convenient than formula (\ref{formr}).
\begin{proposition}
\label{p8}
The kernel $R$ satisfies the recurrence equalities
\begin{gather}
 R^m=\frac{2m-3}{m}\,\b R^{m-1}-\frac{m-3}{m}\,\a R^{m-2}\ \
 \mbox{for}\ \ m\geq 2,
\label{recr1}\\[1mm]
 R^0=-2a_1\sqrt{\de},\ \
 R^1=2a_1\sqrt{\de}\b.
\label{recr2}
\end{gather}
\end{proposition}
\par {\bf Proof}.
The form of formula (\ref{formr1}) differs from the corresponding one in \cite{DZZ09}, Proposition 5.7 only by a constant multiplier.
Therefore Proposition 5.8 in \cite{DZZ09} implies a recurrence formula
\[
 R^m=\frac{2m-3}{m}\,\varkappa\mu R^{m-1}-\frac{m-3}{m}\,\varkappa^2 R^{m-2},\ \ m\geq 2.
\]
Inserting $\vk\mu=\b$ and $\vk^2=\a$ leads to (\ref{recr1}).
Formulas (\ref{recr2}) straightforwardly follow from (\ref{formr}) and (\ref{pmvk}), (\ref{pmvk0}) for $m=1$.
\par
\begin{figure}[h!]
\centering
\includegraphics[height=5cm, width=6cm]{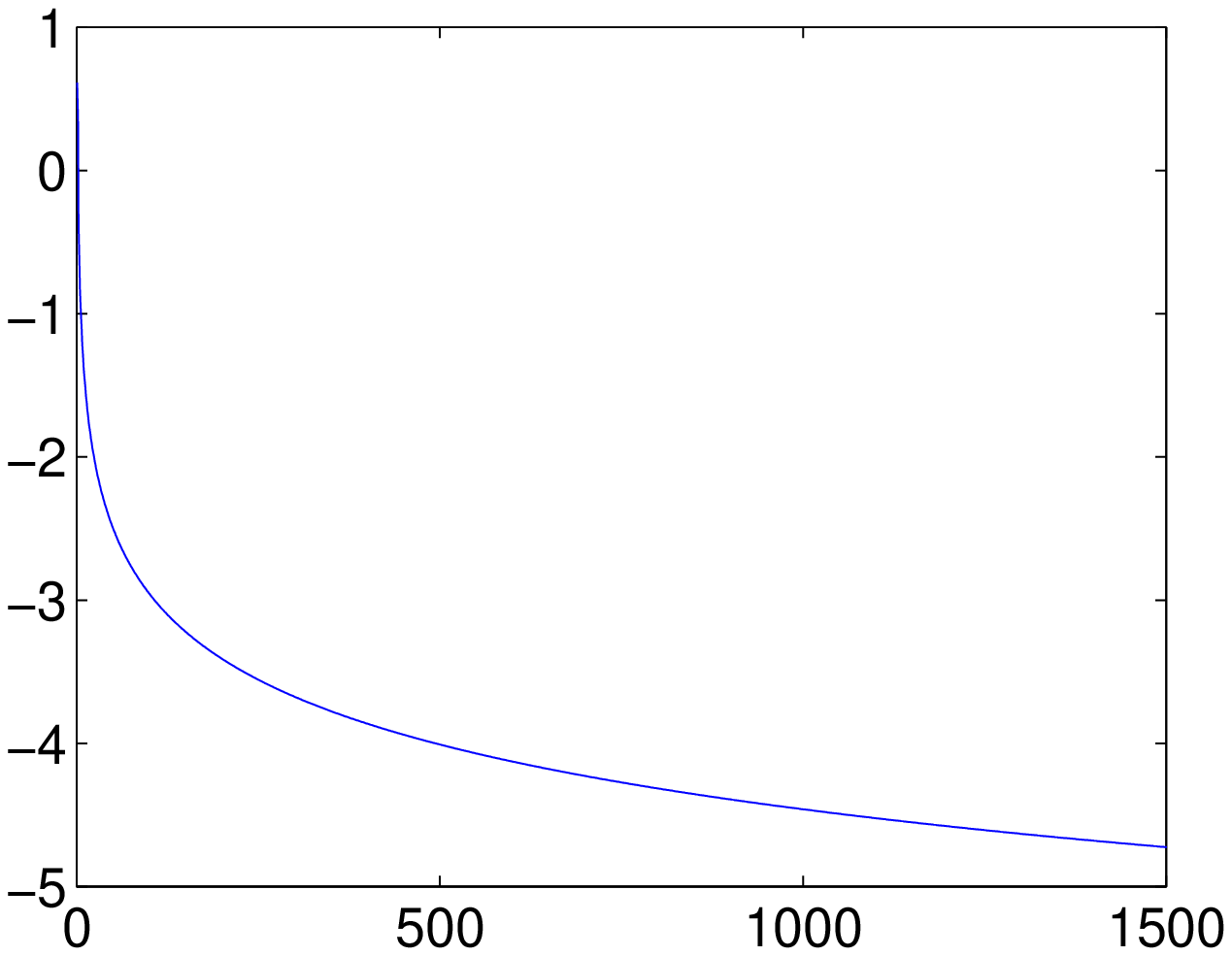}
\includegraphics[height=5cm, width=6cm]{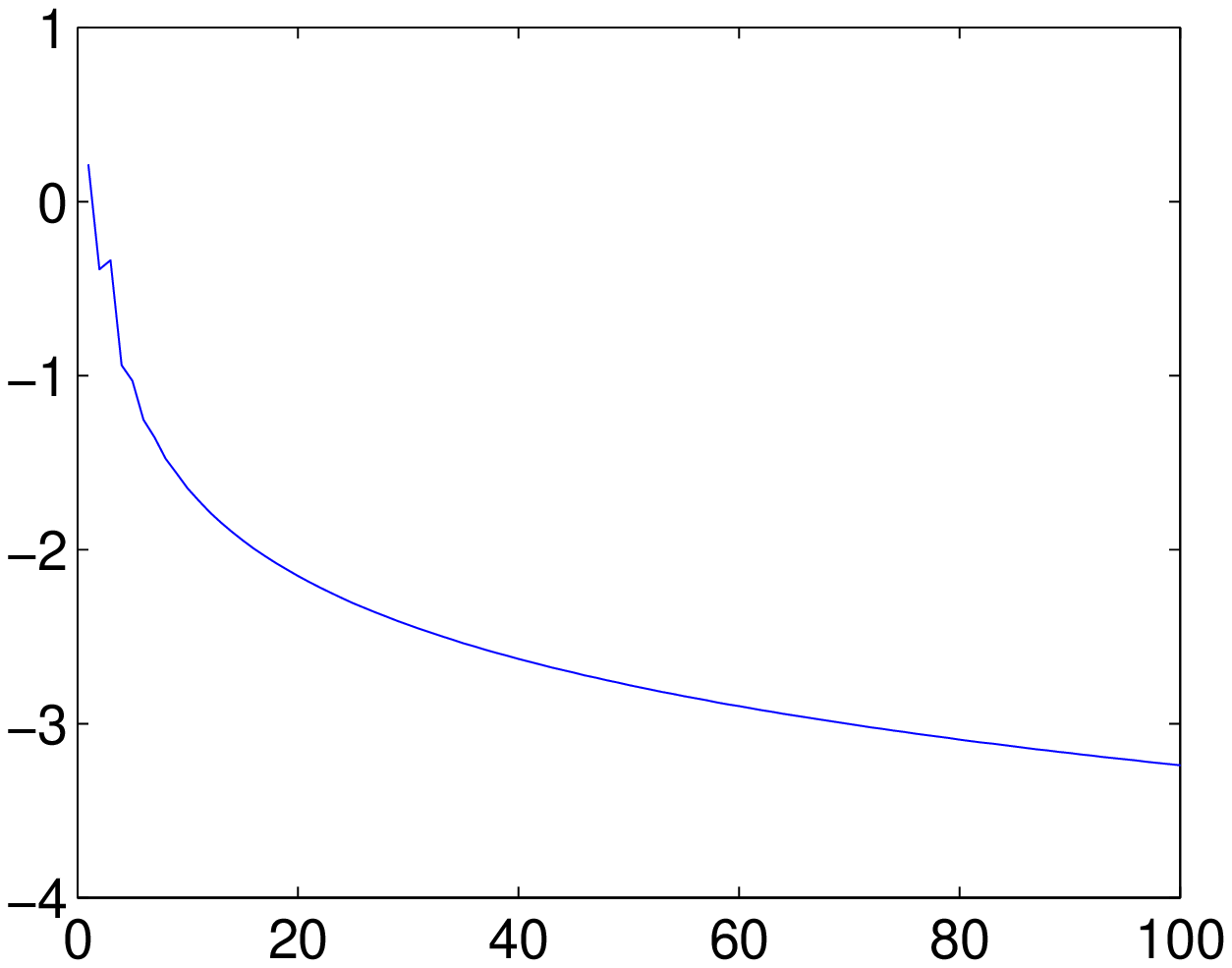}
\caption{Graphs of $\lg|R^m|$ in dependence with $m$, for $\theta=1/12$, in Example 1 (left) and Example 2 (right)}
\label{fig78}
\end{figure}

 Typical graphs of $\lg|R^m|$ are presented on Figures \ref{fig78} for Examples 1 and 2, see Section \ref{sect6} below (the values of parameters are given on Figures \ref{fig12} and \ref{fig56}).
\par One can rather easily extend the above results are to the case of the third boundary condition at $x=0$, or to the Cauchy problem where equation (\ref{sch}) is posed on $\mathbb{R}$. The case $\sigma<1/2$ could be also analyzed under suitable additional condition between steps $\tau$ and $h_j$.
\section{Numerical experiments}
\label{sect6}
\begin{figure}[h]
\centering
\includegraphics[height=5cm, width=6cm]{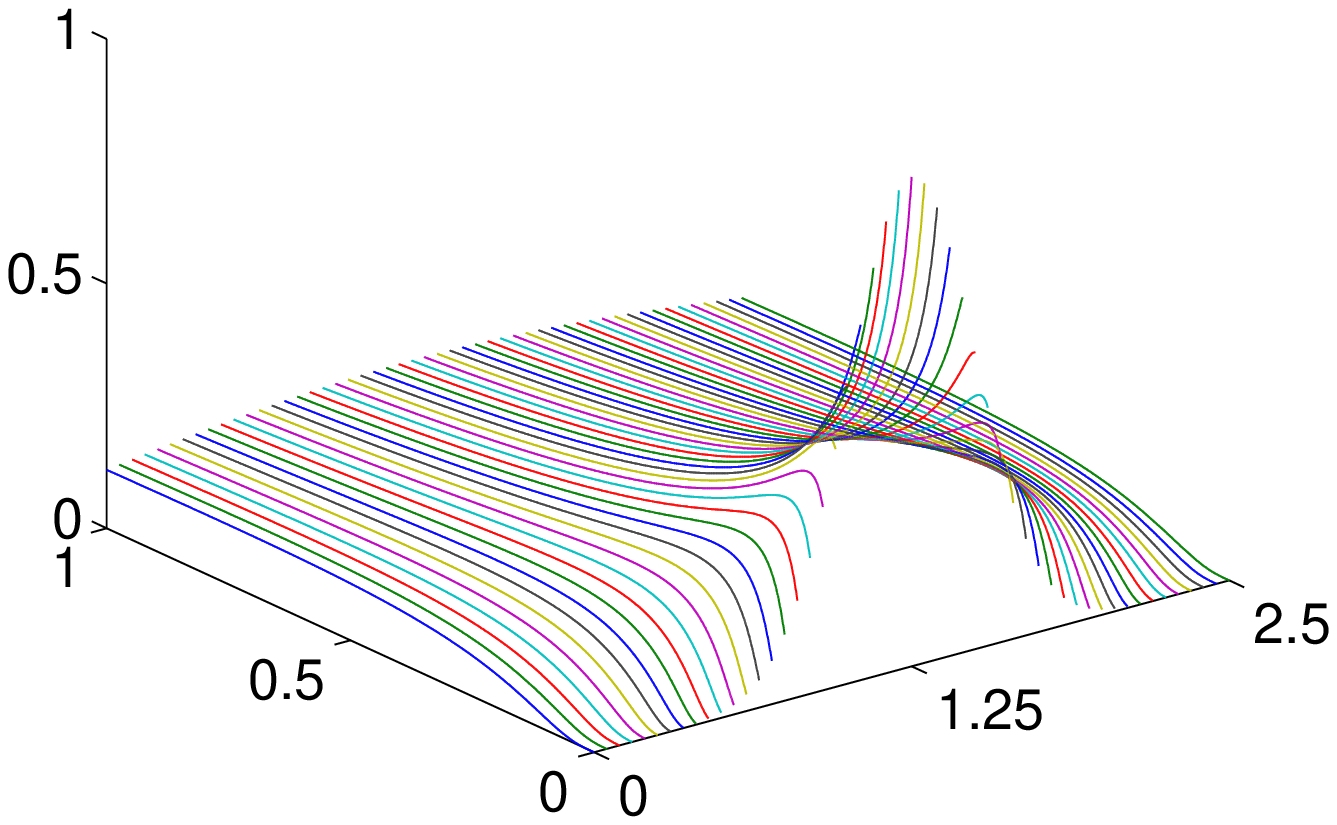}
\includegraphics[height=5cm, width=6cm]{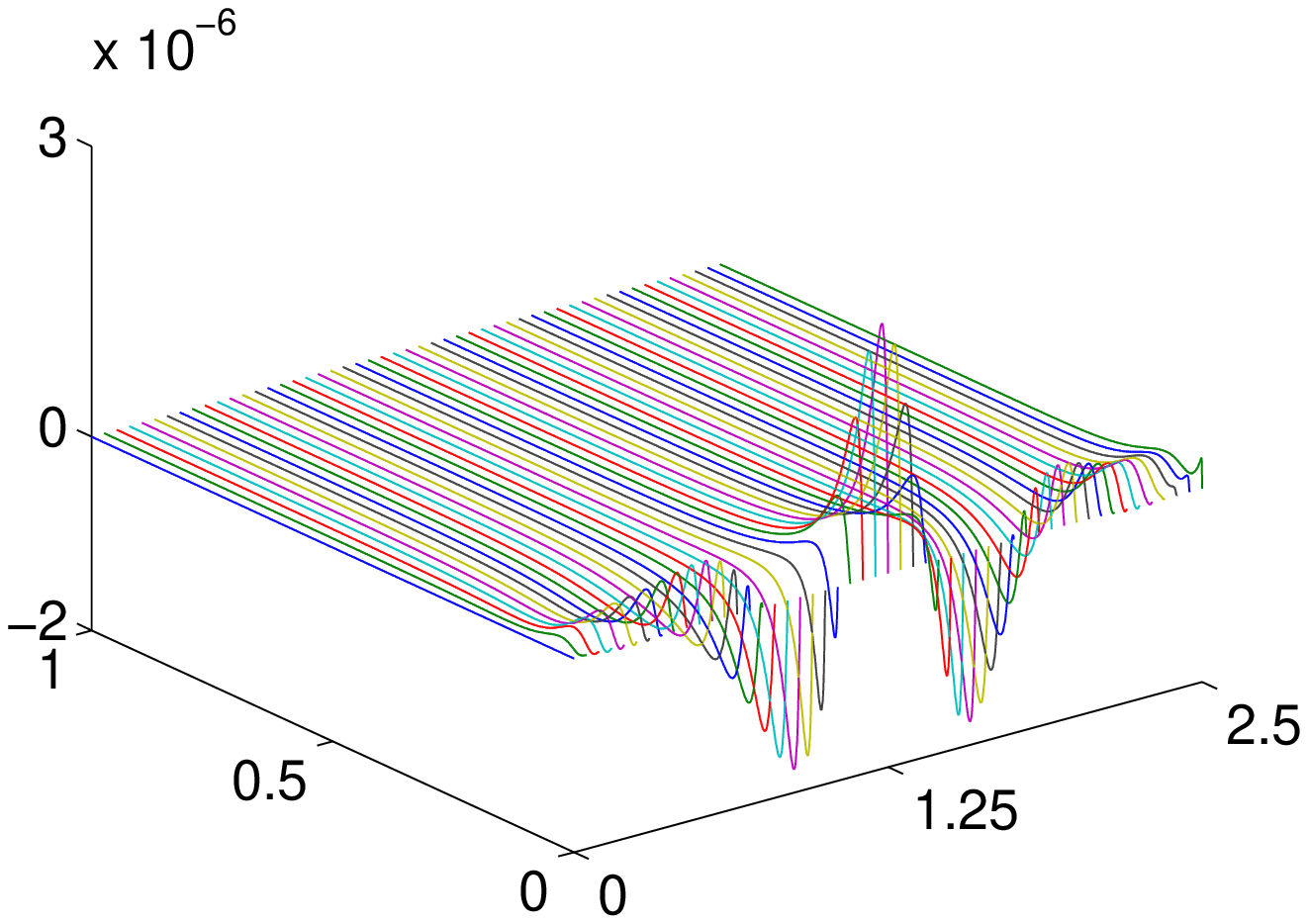}
\caption{Example 1. The numerical solution (left) and its error (right)  for $\sigma=\frac12$, $\theta=\frac{1}{12}$, $h=0.05$ and $\tau=\frac{1}{1500}$}
\label{fig12}
\end{figure}

Consider the initial-boundary value problem (\ref{sch})-(\ref{ic})
for the simplest homogeneous heat equation where $\rho(x)\equiv 1$, $b(x)\equiv 1$, $c(x)\equiv 0$ and $f(x,t)\equiv 0$, for $0\leq t\leq T$.
In Example 1, we base upon an exact solution
\[
 u_1(x,t)=\sqrt{\frac{t_0}{t+t_0}}\,e^{-(x-x_*)^2/[4(t+t_0)]}
\]
with parameters $x_*>0$ and $t_0>0$ and
take the data $g(t)=u_1(0,t)$ and $u^0(x)=u_1(x,0)$.
We choose $x_*=1.25$, $t_0=0.03125$, $T=1$ and $X=2.5$.
Note that $|u^0(x)|< 3.8\cdot 10^{-6}$ for $x\geq X$.
On Fig. \ref{fig12}, we demonstrate the numerical solution and its error computed for $\sigma=\frac12$, $\theta=\frac{1}{12}$, $h=0.05$ and $\tau=\frac{1}{1500}$.
\par In Table \ref{tab1}, the absolute errors are given in dependence with $M$ (where $M\tau=T$) for $\theta=0$ and $\theta=\frac{1}{12}$. For $\theta=0$, they do not practically change for $M\geq 500$ whereas for $\theta=\frac{1}{12}$ the error continues to decrease up to $M=2000$ thus allowing to reach values of $725$ times less.
Note that actually the value $U_J^0\neq 0$ has been used and the minimal reached absolute error is less than this one.
\begin{table}[ht]
\begin{center}
\begin{tabular}{|c|c|c|c|c|c|c|c|}
  \hline
  $M$ & 20 & 50 & 100 & 200 & 500 & 1000 & 2000\\
  \hline
  $\theta=0$ & $4.39\cdot10^{-2}$ & $8.30\cdot10^{-3}$ & $1.30\cdot10^{-3}$ & $5.38\cdot10^{-4}$ & $8.63\cdot10^{-4}$ & $9.15\cdot10^{-4}$ & $9.28\cdot 10^{-4}$\\
  \hline
  $\theta=\frac{1}{12}$ & $4.53\cdot10^{-2}$ & $9.50\cdot10^{-3}$ & $2.20\cdot10^{-3}$ & $4.89\cdot10^{-4}$ & $7.28\cdot10^{-5}$ & $1.35\cdot10^{-5}$ & $1.28\cdot 10^{-6}$\\
  \hline
\end{tabular}
\end{center}
\caption{Example 1. The absolute errors in dependence with $M$ for $h=0.05$}
\label{tab1}
\end{table}
\par Notice that if one sets simply the Neumann boundary condition (i.e., takes $\mathcal{S}=0$) instead of the discrete TBC at $x=X$, then it is necessary to increase $X$ three times to reach the error of the same order of smallness (for the same $h$ and $\tau$); herewith the maximum absolute error is reached at $(x,t)=(X,T)$ (the corresponding graphs are omitted).
\par In Example 2, we take the data $g(t)=t^2$ and $u^0(x)=0$. The exact solution to such problem $u_2(x,t)=32t^2 I_4(\xi)$ is also known and is calculated by applying the recurrence formulas \cite{KE}
\begin{gather*}
I_0(\xi)=\erfc(\xi)
 \equiv\frac{2}{\sqrt{\pi}}\int_{\xi}^{\infty} e^{-\zeta^2}d\zeta,\ \
I_1(\xi)
=\frac{1}{\sqrt{\pi}}e^{-\xi^2}
-\xi\erfc(\xi),
\\[1mm]
I_n(\xi)=\frac{1}{2n}I_{n-2}(\xi)
-\frac{\xi}{n}I_{n-1}(\xi),\ \ n=2,3,4,
\end{gather*}
where $\xi=\frac{x}{2\sqrt{t}}$. Choose $T=1$ and $X=1$.
\par
\begin{figure}[h!]
\centering
\includegraphics[height=5cm, width=6cm]{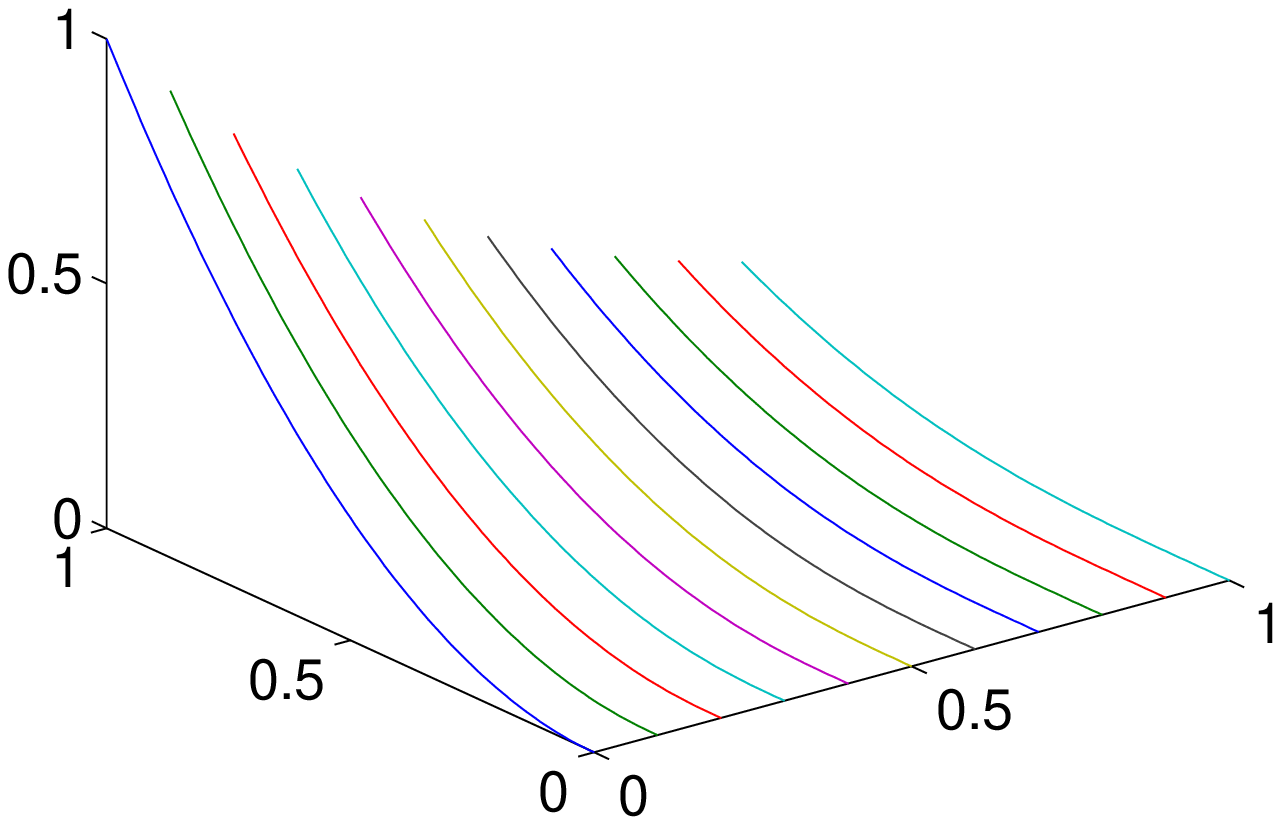}
\includegraphics[height=5cm, width=6cm]{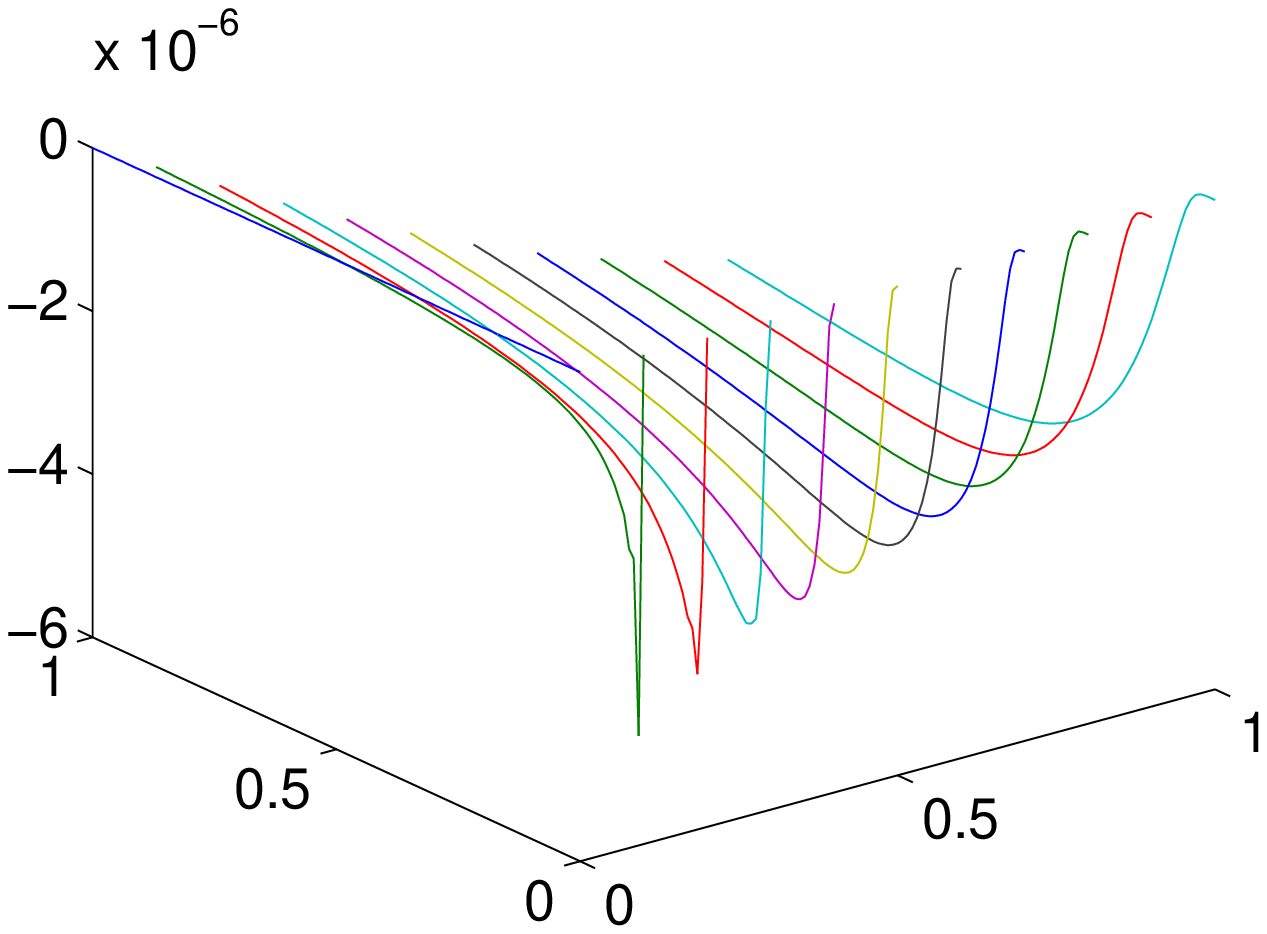}
\caption{Example 2. The numerical solution (left) and its error (right)  for $\sigma=\frac12$, $\theta=\frac{1}{12}$, $h=0.1$ and $\tau=0.01$}
\label{fig56}
\end{figure}

 On Fig. \ref{fig56} we present the numerical solution and its error computed for $\sigma=\frac12$, $\theta=\frac{1}{12}$, $h=0.1$ and $\tau=0.01$. In contrast to Example 1, now the solution is not close to $0$ for $x=X$. The error is maximal
at the node $(x_j,t_m)=(h,\tau)$ (but not on the artificial boundary $x=X$).
Notice that decreasing of $X$ down to $0.2$ for the same mesh steps does not increase the error (for $u^0=0$ and while applying namely the discrete TBC, this is natural and clear from above).
Moreover, if once again one sets simply the Neumann boundary condition instead of the discrete TBC at $x=X$, then (for the same $h$ and $\tau$) the absolute error equals $0.15$ and is unacceptably large. It decreases to the same values as on Fig. \ref{fig56} only when $X$ increases five times.
\par
\begin{figure}[h!]
\centering
\includegraphics[height=5cm, width=6cm]{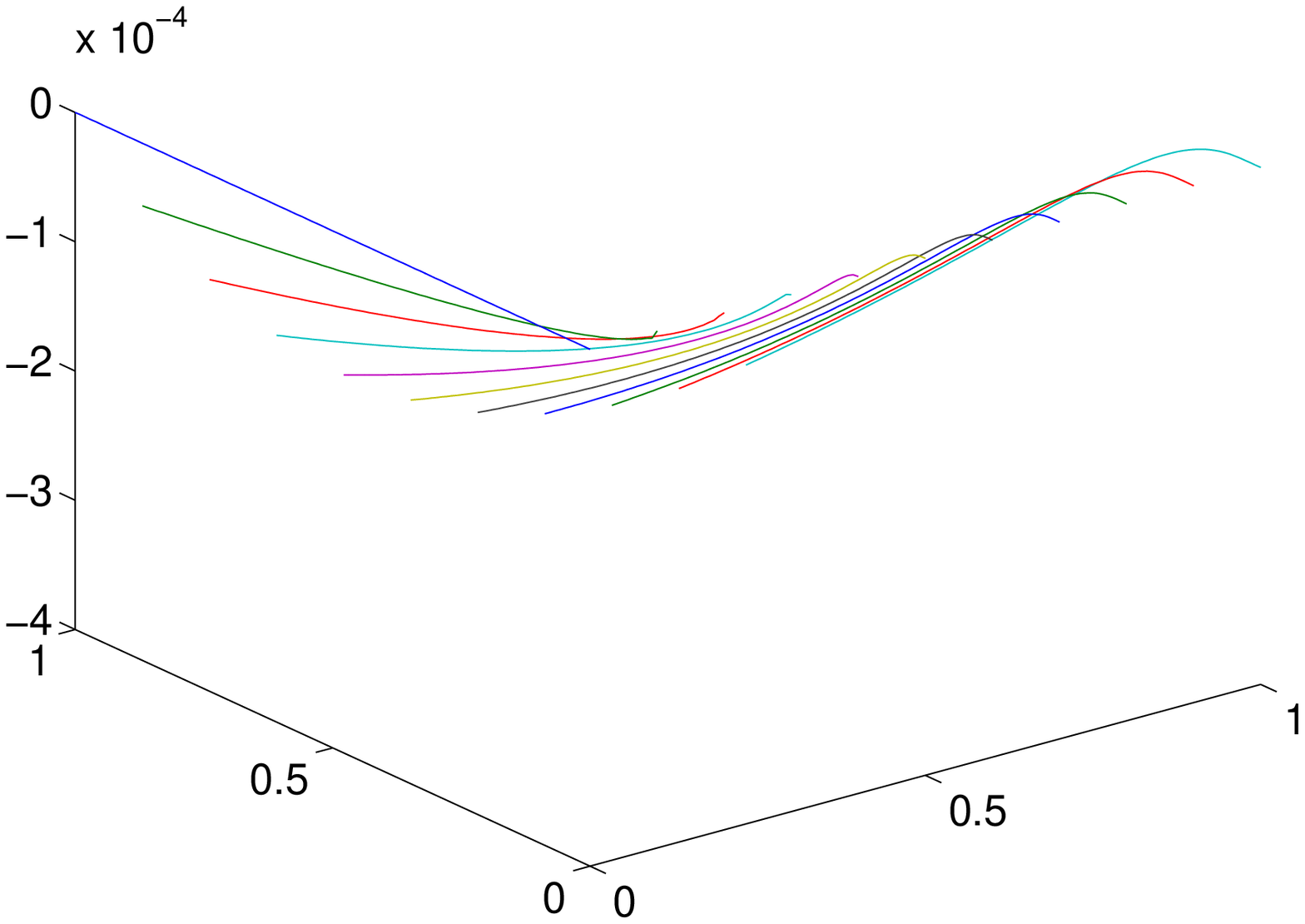}
\includegraphics[height=5cm, width=6cm]{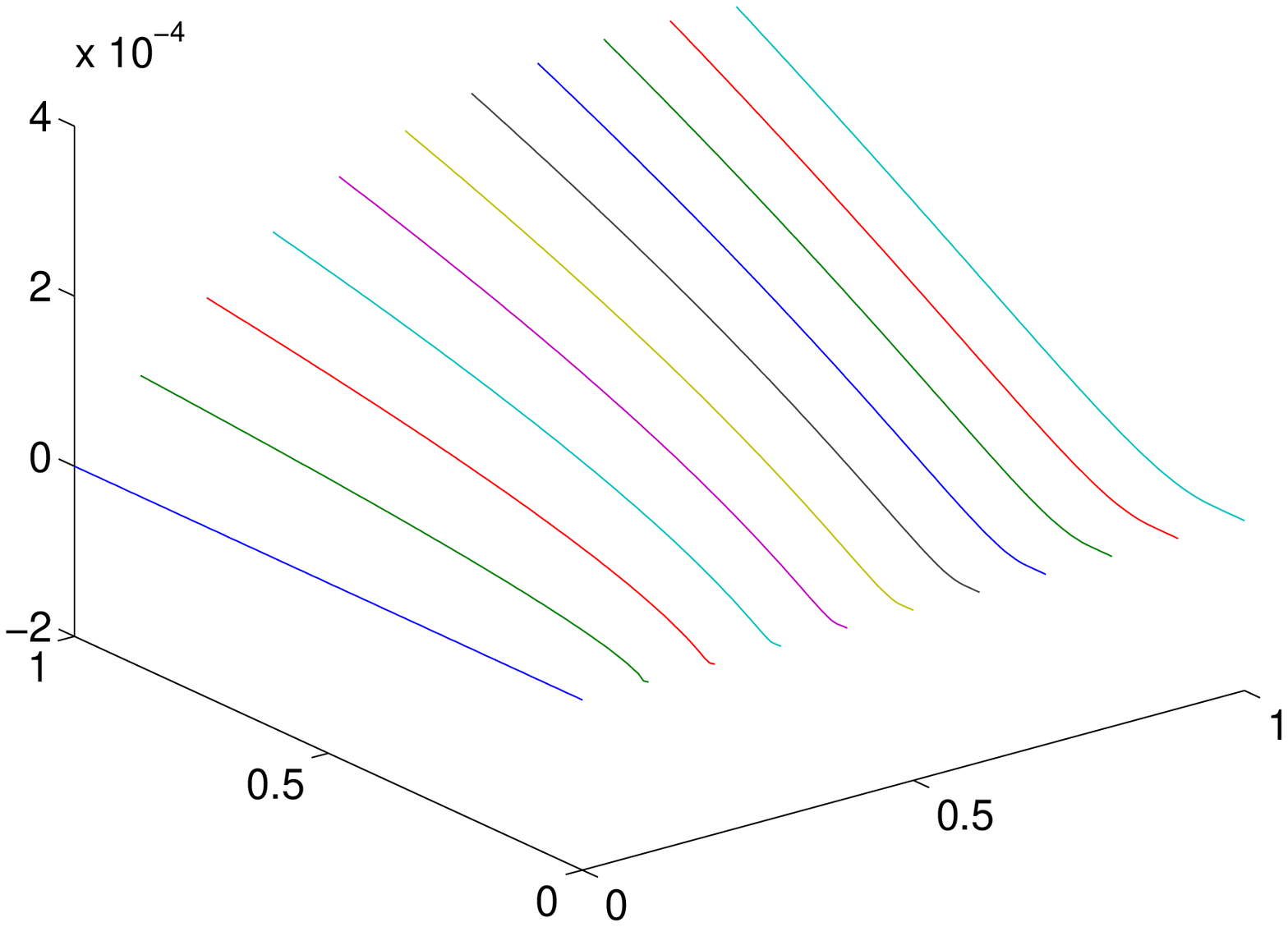}
\caption{Example 2. The error of numerical solution in the cases $\theta=0$ (left) and $\theta=\frac16$ (right) for $\sigma=\frac12$, $h=0.1$ and $\tau=0.01$}
\label{fig56b}
\end{figure}

 On Fig. \ref{fig56b} we give the graphs of errors in the cases $\theta=0$ and $\theta=\frac16$. Their forms are different and the maximum absolute errors are about two orders of magnitude greater than in the case $\theta=\frac{1}{12}$.
\par In addition in Tables \ref{tab2} and \ref{tab3} we put the absolute errors in dependence with $M$ for $h=0.1$ and various $\theta$.
For $\theta=0,\frac16$ and $\frac14$ the errors do not practically change already for $M\geq 50$; herewith their minimum values for $\theta=0$ and $\theta=\frac16$ are close whereas for $\theta=\frac14$ one is approximately twice larger.
In contrast, for $\theta=\frac{1}{12}$ the error continues to decrease rather rapidly up to $M=500$ thus allowing to reach values about four orders of magnitude less (though it increases slightly when $M$ grows further).
\begin{table}[ht]
\begin{center}
\begin{tabular}{|c|c|c|c|c|c|c|}
  \hline
  $M$ & 5 & 10 & 20 & 50 & 100 & 200 \\
  \hline
  $\theta=0$ & $0.020$ & $6.989\cdot10^{-4}$ & $4.239\cdot10^{-4}$ & $3.506\cdot10^{-4}$ & $3.402\cdot10^{-4}$ & $3.376\cdot10^{-4}$ \\
  \hline
  $\theta=\frac{1}{12}$ & $0.020$ & $4.948\cdot10^{-4}$ & $1.224\cdot10^{-4}$ & $1.799\cdot10^{-5}$ & $4.700\cdot10^{-6}$ & $9.300\cdot10^{-7}$\\
  \hline
  $\theta=\frac{1}{6}$ & $0.019$ & $4.654\cdot10^{-4}$ & $2.565\cdot10^{-4}$ & $3.237\cdot10^{-4}$ & $3.340\cdot10^{-4}$ & $3.366\cdot10^{-4}$ \\
\hline
$\theta=\frac{1}{4}$ & $0.019$ & $4.355\cdot10^{-4}$ & $5.931\cdot10^{-4}$ & $6.613\cdot10^{-4}$ & $6.717\cdot10^{-4}$ & $6.743\cdot10^{-4}$\\
\hline
\end{tabular}
\end{center}
\caption{Example 2. The absolute errors in dependence with $M$ for $h=0.1$}
\label{tab2}
\end{table}
\begin{table}[ht]
\begin{center}
\begin{tabular}{|c|c|c|c|c|c|c|c|c|}
\hline
$M$ & 300 & 400 & 500 & 600 & 650 \\
\hline
$\theta=\frac{1}{12}$ & $2.638\cdot10^{-7}$ & $6.797\cdot10^{-8}$ & $4.032\cdot10^{-8}$ & $8.962\cdot10^{-8}$ & $1.062\cdot10^{-7}$\\
\hline
\end{tabular}
\end{center}
\caption{Example 2. The absolute errors in dependence with $M\geq 300$ for $h=0.1$ and $\theta=\frac{1}{12}$}
\label{tab3}
\end{table}
\par The study is carried out by the first author within The National Research University Higher School of Economics' Academic Fund Program in 2012-2013, research grants No. 11-01-0051.
Both authors are also supported by the Federal Agency for Science and Innovations (state contract 14.740.11.0875) and the Russian Foundation for Basic Research, project 10-01-00136.
\makeatletter
\renewcommand{\@biblabel}[1]{#1.\hfill}

\end{document}